\documentclass{article}

\title{\LARGE \textbf{A Sharp Bound for the Circumference in $t$-tough graphs with $t>1$}}
\author{Zh.G. Nikoghosyan\footnote{G.G. Nicoghossian (up to 1997)}}

\begin{document}

\maketitle

\begin{abstract}
It is proved that if  $G$ is a $t$-tough graph of order $n$ and  minimum degree $\delta$ with $t>1$ then either $G$ has a cycle of length at least $\min\{n,2\delta+5\}$ or $G$ is the Petersen graph.  \\


\end{abstract}

\section{Introduction}

Only finite undirected graphs without loops or multiple edges are considered. We reserve $n$, $\delta$, $\kappa$, $c$ and $\tau$ to denote the number of vertices (order),  the minimum degree, connectivity, circumference and the toughness of a graph, respectively. A good reference for any undefined terms is \cite{[2]}. 

The earliest lower bound for the circumference  was developed in 1952 due to Dirac \cite{[3]}.  \\

\noindent\textbf{Theorem A \cite{[3]}.} In every 2-connected graph, $c\ge \min\{n,2\delta\}$. \\

In 1986, Bauer and Schmeichel \cite{[1]} proved that the bound $2\delta$ in Theorem A can be enlarged to $2\delta+2$ by replacing the 2-connectivity condition with 1-toughness.\\

\noindent\textbf{Theorem B \cite{[1]}.} In every 1-tough graph, $c\ge \min\{n,2\delta+2\}$.   \\

In this paper we prove that in Theorem B the bound $2\delta+2$ itself  can be enlarged up to $2\delta+5$ if $\tau>1$ and $G$ is not the Petersen graph.  \\

\noindent\textbf{Theorem 1}. Let $G$ be a  graph with $\tau > 1$. Then either $c\ge\min\{n,2\delta+5\}$ or  $G$ is the Petersen graph.\\

The next result follows immediately.\\

\noindent\textbf{Corollary 1}. Let $G$ be a  graph with $\tau > 1$. If $\delta\ge(n-5)/2$ then $G$ either is hamiltonian or is the Petersen graph. \\

To prove Theorem 1, we need the following result due to Voss \cite{[4]}.  \\

\noindent\textbf{Theorem C \cite{[4]}}. Let $G$ be a hamiltonian graph, $\{v_1,v_2,...,v_t\}\subseteq V(G)$ and $d(v_i)\ge t$ $(i=1,2,...,t)$. Then each pair $x,y$ of vertices of $G$ is connected in $G$ by a path of length at least $t$.\\

\section{Notations and preliminaries}

The set of vertices of a graph $G$ is denoted by $V(G)$ and the set of edges by $E(G)$. For $S$ a subset of $V(G)$, we denote by $G\backslash S$ the maximum subgraph of $G$ with vertex set $V(G)\backslash S$. We write $G[S]$ for the subgraph of $G$ induced by $S$. For a subgraph $H$ of $G$ we use $G\backslash H$ short for $G\backslash V(H)$. The neighborhood of a vertex $x\in V(G)$ will be denoted by $N(x)$.  Furthermore, for a subgraph $H$ of $G$ and $x\in V(G)$, we define $N_H(x)=N(x)\cap V(H)$ and  $d_H(x)=|N_H(x)|$. Let $s(G)$ denote the number of components of a graph $G$. A graph $G$ is $t$-tough if $|S|\ge ts(G\backslash S)$ for every subset $S$ of the vertex set $V(G)$ with $s(G\backslash S)>1$. The toughness of $G$, denoted $\tau(G)$, is the maximum value of $t$ for which $G$ is $t$-tough (taking $\tau(K_n)=\infty$ for all $n\ge 1$).

A simple cycle (or just a cycle) $C$ of length $t$ is a sequence $v_1v_2...v_tv_1$ of distinct vertices $v_1,...,v_t$ with $v_iv_{i+1}\in E(G)$ for each $i\in \{1,...,t\}$, where $v_{t+1}=v_1$. When $t=2$, the cycle $C=v_1v_2v_1$ on two vertices $v_1, v_2$ coincides with the edge $v_1v_2$, and when $t=1$, the cycle $C=v_1$ coincides with the vertex $v_1$. So, all vertices and edges in a graph can be considered as cycles of lengths 1 and 2, respectively. A graph $G$ is hamiltonian if $G$ contains a Hamilton cycle, i.e. a cycle of length $n$. A cycle $C$ in $G$ is dominating if  $G\backslash C$ is edgeless.

Paths and cycles in a graph $G$ are considered as subgraphs of $G$. If $Q$ is a path or a cycle, then the length of $Q$, denoted by $|Q|$, is $|E(Q)|$. We write $Q$ with a given orientation by $\overrightarrow{Q}$. For $x,y\in V(Q)$, we denote by $x\overrightarrow{Q}y$ the subpath of $Q$ in the chosen direction from $x$ to $y$. For $x\in V(C)$, we denote the $h$-th successor and the $h$-th predecessor of $x$ on $\overrightarrow{C}$ by $x^{+h}$ and $x^{-h}$, respectively. We abbreviate $x^{+1}$ and $x^{-1}$ by $x^+$ and $x^-$, respectively.  For each $X\subset V(C)$, we define $X^{+h}=\{x^{+h}|x\in X\}$ and $X^{-h}=\{x^{-h}|x\in X\}$. \\

\noindent\textbf{Special definitions}. Let $G$ be a graph, $C$ a longest cycle in $G$ and $P=x\overrightarrow{P}y$ a longest path in $G\backslash C$ of length $\overline{p}\ge0$. Let $\xi_1,\xi_2,...,\xi_s$ be the elements of $N_C(x)\cup N_C(y)$ occuring on $C$ in a consecutive order. Set 
$$
I_i=\xi_i\overrightarrow{C}\xi_{i+1},       \     I_i^\ast=\xi_i^+\overrightarrow{C}\xi_{i+1}^-   \   \  (i=1,2,...,s),
$$
where $\xi_{s+1}=\xi_1$. 

$(1)$ The segments  $I_1,I_2,...,I_s$ are called elementary segments on $C$ induced by $N_C(x)\cup N_C(y)$.

$(2)$ We call a path $L=z\overrightarrow{L}w$ an intermediate path between two distinct elementary segments $I_a$ and $I_b$ if
$$
z\in V(I_a^\ast),   \    w\in V(I_b^\ast),    \    V(L)\cap V(C\cup P)=\{z,w\}.
$$

$(3)$ Define $\Upsilon(I_{i_1},I_{i_2},...,I_{i_t})$ to be the set of all intermediate paths between elementary segments  $I_{i_1},I_{i_2},...,I_{i_t}$.\\

\noindent\textbf{Lemma 1.} Let $G$ be a graph, $C$ a longest cycle in $G$ and $P=x\overrightarrow{P}y$ a longest path in $G\backslash C$ of length $\overline{p}\ge1$. If  $|N_C(x)|\ge2$, $|N_C(y)|\ge2$ and  $N_C(x)\not=N_C(y)$ then
$$
|C|\ge\left\{ 
\begin{array}{lll}
3\delta+\max\{\sigma_1, \sigma_2\}-1\ge3\delta & \mbox{if} & \mbox{ }\overline{p}=1, \\ \max\{2\overline{p}+8, 4\delta-2\overline{p}\} & \mbox{if} & \mbox{ }%
\overline{p}\ge2, 
\end{array}
\right. 
$$
where $\sigma_1=|N_C(x)\backslash N_C(y)|$ and $\sigma_2=|N_C(y)\backslash N_C(x)|$.\\

\noindent\textbf{Lemma 2.} Let $G$ be a graph, $C$ a longest cycle in $G$ and $P=x\overrightarrow{P}y$ a longest path in $G\backslash C$ of length $\overline{p}\ge0$.  If $N_C(x)=N_C(y)$ and $|N_C(x)|\ge2$ then for each  elementary segments $I_a$ and $I_b$ induced by $N_C(x)\cup N_C(y)$,\\

(a1) if $L$ is an intermediate path between  $I_a$ and $I_b$ then 
$$
|I_a|+|I_b|\ge2\overline{p}+2|L|+4,
$$

(a2) if $\Upsilon(I_a,I_b)\subseteq E(G)$ and $|\Upsilon(I_a,I_b)|=i$\ \  for some $i\in\{1,2,3\}$ then 
$$
|I_a|+|I_b|\ge2\overline{p}+i+5,
$$

(a3) if $\Upsilon(I_a,I_b)\subseteq E(G)$ and $\Upsilon(I_a,I_b)$ contains two independent intermediate edges then 
$$
|I_a|+|I_b|\ge2\overline{p}+8.
$$

\noindent\textbf{Lemma 3}. Let $G$ be a graph and $C$ a longest cycle in $G$. Then either $|C|\ge\kappa(\delta+1)$ or there is a longest path $P=x_1\overrightarrow{P}x_2$ in $G\backslash C$ with $|N_C(x_i)|\ge2$ $(i=1,2)$.

\section{Proofs}

\noindent\textbf{Proof of Lemma 1}. Put
$$
A_1=N_C(x)\backslash N_C(y),    \    A_2=N_C(y)\backslash N_C(x), \   M=N_C(x)\cap N_C(y).
$$
By the hypothesis,  $N_C(x)\not=N_C(y)$, implying that 
$$
\max \{|A_1|,|A_2|\}\ge1.
$$ 
Let $\xi_1,\xi_2,...,\xi_s$ be the elements of $N_C(x)\cup N_C(y)$ occuring on $C$ in a consecutive order. Put $I_i=\xi_i\overrightarrow{C}\xi_{i+1}$ $(i=1,2,...,s)$, where $\xi_{s+1}=\xi_1$. Clearly, $s=|A_1|+|A_2|+|M|$. Since $C$ is extreme, $|I_i|\ge2$ $(i=1,2,...,s)$. Next, if $\{\xi_i,\xi_{i+1}\}\cap M\not=\emptyset$ for some $i\in\{1,2,...,s\}$ then $|I_i|\ge\overline{p}+2$. Further, if either $\xi_i\in A_1$, $\xi_{i+1}\in A_2$ or $\xi_i\in A_2$, $\xi_{i+1}\in A_1$ then again $|I_i|\ge\overline{p}+2$. \\

\textbf{Case 1}. $\overline{p}=1$.

\textbf{Case 1.1}. $|A_i|\ge1$ $(i=1,2)$.

It follows that among  $I_1,I_2,...,I_s$ there are  $|M|+2$ segments of length at least $\overline{p}+2$. Observing also that each of the remaining $s-(|M|+2)$ segments has a length at least 2, we have 
$$
|C|\ge(\overline{p}+2)(|M|+2)+2(s-|M|-2)
$$
$$
=3(|M|+2)+2(|A_1|+|A_2|-2)
$$
$$
=2|A_1|+2|A_2|+3|M|+2.
$$

Since $|A_1|=d(x)-|M|-1$ and  $|A_2|=d(y)-|M|-1$, 
$$
|C|\ge2d(x)+2d(y)-|M|-2\ge3\delta+d(x)-|M|-2.
$$
Recalling that $d(x)=|M|+|A_1|+1$, we get 
$$
|C|\ge3\delta+|A_1|-1=3\delta+\sigma_1-1.
$$
Analogously, $|C|\ge3\delta+\sigma_2-1$. So, 
$$
|C|\ge3\delta+\max \{\sigma_1,\sigma_2\}-1\ge3\delta.
$$

\textbf{Case 1.2}. Either $|A_1|\ge1, |A_2|=0$ or $|A_1|=0, |A_2|\ge1$.

Assume w.l.o.g. that $|A_1|\ge1$ and $|A_2|=0$, i.e. $|N_C(y)|=|M|\ge2$ and $s=|A_1|+|M|$ . Hence, among $I_1,I_2,...,I_s$ there are $|M|+1$ segments of length at least $\overline{p}+2=3$. Taking into account that each of the remaining $s-(|M|+1)$ segments has a length at least 2 and $|M|+1=d(y)$, we get
$$
|C|\ge 3(|M|+1)+2(s-|M|-1)=3d(y)+2(|A_1|-1)
$$
$$
\ge3\delta+|A_1|-1=3\delta+\max\{\sigma_1,\sigma_2\}-1\ge3\delta.
$$

\textbf{Case 2}. $\overline{p}\ge2$.

We first prove that $|C|\ge2\overline{p}+8$. Since $|N_C(x)|\ge2$ and $|N_C(y)|\ge2$,  there are at least two segments among $I_1,I_2,...,I_s$ of length at least $\overline{p}+2$. If $|M|=0$ then clearly $s\ge4$ and 
$$
|C|\ge2(\overline{p}+2)+2(s-2)\ge2\overline{p}+8.
$$
Otherwise, since $\max\{|A_1|,|A_2|\}\ge1$, there are at least three elementary segments of length at least $\overline{p}+2$, that is 
$$
|C|\ge3(\overline{p}+2)\ge2\overline{p}+8.
$$
So, in any case, $|C|\ge 2\overline{p}+8$. 

To prove that $|C|\ge 4\delta-2\overline{p}$, we distinguish two main cases.\\

\textbf{Case 2.1}. $|A_i|\ge1$ $(i=1,2)$.

It follows that among  $I_1,I_2,...,I_s$ there are $|M|+2$ segments of length at least $\overline{p}+2$. Further, since each of the remaining $s-(|M|+2)$ segments has a length at least 2, we get 
$$
|C|\ge (\overline{p}+2)(|M|+2)+2(s-|M|-2)
$$
$$
=(\overline{p}-2)|M|+(2\overline{p}+4|M|+4)+2(|A_1|+|A_2|-2)
$$
$$
\ge2|A_1|+2|A_2|+4|M|+2\overline{p}.
$$
Observing also that 
$$
|A_1|+|M|+\overline{p}\ge d(x),  \quad   |A_2|+|M|+\overline{p}\ge d(y),  
$$
we have 

$$
2|A_1|+2|A_2|+4|M|+2\overline{p}
$$
$$
\ge 2d(x)+2d(y)-2\overline{p}\ge4\delta-2\overline{p},
$$
implying that $|C|\ge4\delta-2\overline{p}$.\\

 \textbf{Case 2.2}. Either $|A_1|\ge1, |A_2|=0$ or $|A_1|=0, |A_2|\ge1$.

Assume w.l.o.g. that $|A_1|\ge1$ and $|A_2|=0$, i.e.  $|N_C(y)|=|M|\ge2$ and $s=|A_1|+|M|$. It follows that among  $I_1,I_2,...,I_s$ there are  $|M|+1$ segments of length at least $\overline{p}+2$. Observing also that $|M|+\overline{p}\ge d(y)\ge\delta$, i.e. $2\overline{p}+4|M|\ge 4\delta-2\overline{p}$, we get 
$$
|C|\ge(\overline{p}+2)(|M|+1)\ge(\overline{p}-2)(|M|-1)+2\overline{p}+4|M|
$$
$$
\ge 2\overline{p}+4|M|\ge4\delta-2\overline{p}.       \quad            \quad                    \rule{7pt}{6pt} 
$$

\noindent\textbf{Proof of Lemma 2}. Let $\xi_1,\xi_2,...,\xi_s$ be the elements of $N_C(x)$ occuring on $C$ in a consecutive order. Put $I_i=\xi_i\overrightarrow{C}\xi_{i+1}$ $(i=1,2,...,s)$, where $\xi_{s+1}=\xi_1.$  To prove $(a1)$, let $L=z\overrightarrow{L}w$ be an intermediate path between elementary segments $I_a$ and $I_b$ with $z\in V(I_a^\ast)$ and   $w\in V(I_b^\ast)$.
Put
$$
|\xi_a\overrightarrow{C}z|=d_1,  \   |z\overrightarrow{C}\xi_{a+1}|=d_2,      \      |\xi_b\overrightarrow{C}w|=d_3,       \     |w\overrightarrow{C}\xi_{b+1}|=d_4,
$$
$$
C^\prime=\xi_ax\overrightarrow{P}y\xi_b\overleftarrow{C}z\overrightarrow{L}w\overrightarrow{C}\xi_a.
$$
Clearly, 
$$
|C^\prime|=|C|-d_1-d_3+|L|+|P|+2.           
$$
Since $C$ is extreme, we have $|C|\ge|C^\prime|$, implying that $d_1+d_3\ge\overline{p}+|L|+2$.  By a symmetric argument, $d_2+d_4\ge\overline{p}+|L|+2$. Hence 
$$
|I_a|+|I_b|=\sum_{i=1}^4d_i\ge2\overline{p}+2|L|+4.
$$

The proof of $(a1)$ is complete. To proof  $(a2)$ and $(a3)$, let  $\Upsilon(I_a,I_b)\subseteq E(G)$ and $|\Upsilon(I_a,I_b)|=i$ for some $i\in \{1,2,3\}$.\\

\textbf{Case 1}. $i=1$.

It follows that $\Upsilon(I_a,I_b)$ consists of a unique intermediate edge $L=zw$. By (a1), 
$$
|I_a|+|I_b|\ge2\overline{p}+2|L|+4=2\overline{p}+6.
$$

\textbf{Case 2}. $i=2$.

It follows that $\Upsilon(I_a,I_b)$ consists of two edges $e_1,e_2$. Put $e_1=z_1w_1$ and $e_2=z_2w_2$, where $\{z_1,z_2\}\subseteq V(I_a^\ast)$ and $\{w_1,w_2\}\subseteq V(I_b^\ast)$.\\

\textbf{Case 2.1}. $z_1\not=z_2$ and $w_1\not=w_2$.

Assume w.l.o.g. that $z_1$ and $z_2$ occur in this order on $I_a$. \\

\textbf{Case 2.1.1}. $w_2$ and  $w_1$ occur in this order on $I_b$.

Put
$$
|\xi_a\overrightarrow{C}z_1|=d_1,  \   |z_1\overrightarrow{C}z_2|=d_2,      \    |z_2\overrightarrow{C}\xi_{a+1}|=d_3, 
$$
$$
|\xi_b\overrightarrow{C}w_2|=d_4,       \       |w_2\overrightarrow{C}w_1|=d_5,         \        |w_1\overrightarrow{C}\xi_{b+1}|=d_6, 
$$
$$
C^{\prime}=\xi_a\overrightarrow{C}z_1w_1\overleftarrow{C}w_2z_2\overrightarrow{C}\xi_b x\overrightarrow{P}y\xi_{b+1}\overrightarrow{C}\xi_a.
$$
Clearly, 
$$
|C^{\prime}|=|C|-d_2-d_4-d_6+|\{e_1\}|+|\{e_2\}|+|P|+2
$$
$$
=|C|-d_2-d_4-d_6+\overline{p}+4.
$$
Since $C$ is extreme, $|C|\ge |C^{\prime}|$, implying that $d_2+d_4+d_6\ge \overline{p}+4$. By a symmetric argument, $d_1+d_3+d_5\ge\overline{p}+4$. Hence
$$
|I_a|+|I_b|= \sum_{i=1}^6d_i\ge2\overline{p}+8.
$$

\textbf{Case 2.1.2}. $w_1$ and $w_2$ occur in this order on $I_b$.

Putting
$$
C^{\prime}=\xi_a\overrightarrow{C}z_1w_1\overrightarrow{C}w_2z_2\overrightarrow{C}\xi_b x\overrightarrow{P}y\xi_{b+1}\overrightarrow{C}\xi_a,
$$
we can argue as in Case 2.1.1. \\

\textbf{Case 2.2}. Either $z_1=z_2$, $w_1\not=w_2$ or $z_1\not=z_2$, $w_1=w_2$.

Assume w.l.o.g. that $z_1\not=z_2$, $w_1=w_2$ and $z_1, z_2$ occur  in this order on $I_a$. Put
$$
|\xi_a\overrightarrow{C}z_1|=d_1,  \   |z_1\overrightarrow{C}z_2|=d_2,      \    |z_2\overrightarrow{C}\xi_{a+1}|=d_3, 
$$
$$
|\xi_b\overrightarrow{C}w_1|=d_4,       \           |w_1\overrightarrow{C}\xi_{b+1}|=d_5, 
$$
$$
C^{\prime}=\xi_a x\overrightarrow{P}y\xi_b\overleftarrow{C}z_1w_1\overrightarrow{C}\xi_a,
$$
$$
C^{\prime\prime}=\xi_a\overrightarrow{C}z_2w_1\overleftarrow{C}\xi_{a+1}x\overrightarrow{P}y\xi_{b+1}\overrightarrow{C}\xi_a.
$$

Clearly, 
$$
|C^{\prime}|=|C|-d_1-d_4+|\{e_1\}|+|P|+2=|C|-d_1-d_4+\overline{p}+3,
$$
$$
|C^{\prime\prime}|=|C|-d_3-d_5+|\{e_2\}|+|P|+2=|C|-d_3-d_5+\overline{p}+3.
$$
Since $C$ is extreme, $|C|\ge |C^{\prime}|$ and $|C|\ge |C^{\prime\prime}|$, implying that 
$$
d_1+d_4\ge \overline{p}+3,  \    d_3+d_5\ge \overline{p}+3. 
$$
Hence, 
$$
|I_a|+|I_b|= \sum_{i=1}^5d_i\ge d_1+d_3+d_4+d_5+1\ge2\overline{p}+7.
$$

\textbf{Case 3}. $i=3$.

It follows that $\Upsilon(I_a,I_b)$ consists of three edges $e_1,e_2,e_3$. Let $e_i=z_iw_i$ $(i=1,2,3)$, where $\{z_1,z_2,z_3\}\subseteq V(I_a^\ast)$ and $\{w_1,w_2,w_3\}\subseteq V(I_b^\ast)$. If there are two independent edges among $e_1,e_2,e_3$ then we can argue as in Case 2.1. Otherwise, we can assume w.l.o.g. that $w_1=w_2=w_3$ and $z_1,z_2,z_3$ occur in this order on $I_a$. Put
$$
|\xi_a\overrightarrow{C}z_1|=d_1,  \   |z_1\overrightarrow{C}z_2|=d_2,      \    |z_2\overrightarrow{C}z_3|=d_3,      
$$
$$
|z_3\overrightarrow{C}\xi_{a+1}|=d_4,      \    |\xi_b\overrightarrow{C}w_1|=d_5,       \           |w_1\overrightarrow{C}\xi_{b+1}|=d_6, 
$$
$$
C^{\prime}=\xi_a x\overrightarrow{P}y\xi_b\overleftarrow{C}z_1w_1\overrightarrow{C}\xi_a,
$$
$$
C^{\prime\prime}=\xi_a\overrightarrow{C}z_3w_1\overleftarrow{C}\xi_{a+1}x\overrightarrow{P}y\xi_{b+1}\overrightarrow{C}\xi_a.
$$

Clearly, 
$$
|C^{\prime}|=|C|-d_1-d_5+|\{e_1\}|+\overline{p}+2,
$$
$$
|C^{\prime\prime}|=|C|-d_4-d_6+|\{e_3\}|+\overline{p}+2.
$$
Since $C$ is extreme, we have $|C|\ge |C^{\prime}|$ and $|C|\ge |C^{\prime\prime}|$, implying that 
$$
d_1+d_5\ge \overline{p}+3,   \    d_4+d_6\ge \overline{p}+3. 
$$
Hence, 
$$
|I_a|+|I_b|= \sum_{i=1}^6d_i\ge d_1+d_4+d_5+d_6+2\ge2\overline{p}+8.     \quad      \quad    \rule{7pt}{6pt} 
$$

\noindent\textbf{Proof of Lemma 3}. Choose a longest path $P=x_1\overrightarrow{P}x_2$ in $G\backslash C$ so as to maximize $|N_C(x_1)|$. Let $y_1,...,y_t$ be the elements of $N^+_P(x_2)$ occuring on $P$ in a consecutive order. Put
$$
P_i=x_1\overrightarrow{P}y_i^-x_2\overleftarrow{P}y_i  \  (i=1,...,t),  \  \  H=G[V(y_1^-\overrightarrow{P}x_2)].
$$ 
Since $P_i$ is a longest path in $G\backslash C$ for each $i\in\{1,...,t\}$, we can assume w.l.o.g. that $P$ is chosen such that $|V(H)|$ is maximum. It follows in particular that $N_P(y_i)\subseteq V(H)$ $(i=1,...,t)$. \\

\textbf{Case 1}. $|N_C(x_1)|=0$.

Since $|N_C(x_1)|$ is maximum, we have $|N_C(y_i)|=0$ $(i=1,...,t)$, implying that $N(y_i)\subseteq V(H)$ and $d_H(y_i)=d(y_i)\ge\delta$ $(i=1,...,t)$. Further, since $y_t=x_2$, we have $d_P(x_2)\ge\delta$, that is $t\ge\delta$. By Theorem C, for each distinct $u,v\in V(H)$, there is a path in $H$ of length at least $\delta$, connecting $u$ and $v$. Since $H$ and $C$ are connected by at least $\kappa$ vertex disjoint paths, we have $|C|\ge\kappa(\delta+2)$.\\

\textbf{Case 2}. $|N_C(x_1)|=1$.

Since $|N_C(x_1)|$ is maximum, we have $|N_C(y_i)|\le1$ $(i=1,...,t)$, implying that $|N_H(y_i)|\ge\delta-1$ $(i=1,...,t)$, where $t\ge\delta-1$. By Theorem C, $|C|\ge\kappa(\delta+1)$.\\

\textbf{Case 3}. $|N_C(x_1)|\ge2$.

If $|N_C(y_i)|\ge2$ for some $i\in \{1,...,t\}$ then we are done. Otherwise $|N_C(y_i)|\le1$ $(i=1,...,t)$ and, as in Case 2, $|C|\ge\kappa(\delta+1)$.           \quad      \quad    \rule{7pt}{6pt} \\

\noindent\textbf{Proof of Theorem 1}.  If $\kappa\le2$ then clearly $\tau\le1$, contradicting the hypothesis.  Next, if $c\ge2\delta+5$ then we are done. So, we can assume that
$$
\delta\ge\kappa\ge3,   \  \  \  c\le2\delta+4.                                 \eqno{(1)}
$$

Let $C$ be a longest cycle in $G$ and $P=x_1\overrightarrow{P}x_2$ a longest path in $G\backslash C$ of length $\overline{p}$. If $|V(P)|\le 0$ then $C$ is a Hamilton cycle and we are done. Let $|V(P)|\ge1$.  Put $X=N_C(x_1)\cup N_C(x_2)$ and let $\xi_1,...,\xi_s$ be the elements of $X$ occuring on $C$ in a consecutive order. Put
$$
I_i=\xi_i\overrightarrow{C}\xi_{i+1},       \     I_i^\ast=\xi_i^+\overrightarrow{C}\xi_{i+1}^-  \  \ (i=1,...,s),
$$
where $\xi_{s+1}=\xi_1.$   \\

\textbf{Claim 1}. Let $N_C(x_1)=N_C(x_2)$ and let $\xi_a,\xi_b$ be two distinct elements of $X$. If either $|\xi_a\overrightarrow{C}y|+|\xi_b\overrightarrow{C}z|\le\overline{p}+2$ or $|y\overrightarrow{C}\xi_{a+1}|+|z\overrightarrow{C}\xi_{b+1}|\le\overline{p}+2$ for some $y\in V(I_a^\ast)$ and $z\in V(I_b^\ast)$, then $yz\not\in E(G)$.

\textbf{Proof}. Assume the contrary, that is $yz\in E(G)$. If  $|\xi_a\overrightarrow{C}y|+|\xi_b\overrightarrow{C}z|\le\overline{p}+2$ then
$$
|\xi_ax_1\overrightarrow{P}x_2\xi_b\overleftarrow{C}yz\overrightarrow{C}\xi_a|=|C|-|\xi_a\overrightarrow{C}y|-|\xi_b\overrightarrow{C}z|+\overline{p}+3\ge|C|+1,
$$
a contradiction. By a symmetric argument, we reach a contradiction when $|y\overrightarrow{C}\xi_{a+1}|+|z\overrightarrow{C}\xi_{b+1}|\le\overline{p}+2$.   \  \  \  $\Delta$\\

\textbf{Case 1}. $\overline{p}=0$.

It follows that $P=x_1$ and $s=d(x_1)\ge\delta\ge3$. The next claim can be derived from (1) and Lemma 2 easily.\\

\textbf{Claim 2}. $(1)$ $|I_i|+|I_j|\le8$ for each distinct $i,j\in \{1,...,s\}$.

$(2)$ If $|I_a|+|I_b|=8$ for some distinct $a,b\in \{1,...,s\}$ then  $|I_i|=2$ for each $i\in \{1,...,s\}\backslash \{a,b\}$. 

$(3)$ If $|I_a|=6$ for some $a\in\{1,...,s\}$ then $|I_i|=2$ for each $i\in \{1,...,s\}\backslash \{a\}$.  

$(4)$ There are at most four segments of length at least 3.

$(5)$ If $|I_a|\ge3$, $|I_b|\ge3$, $|I_f|\ge3$, $|I_g|\ge3$ for some distinct $a,b,f,g\in\{1,...,s\}$ then $|I_a|=|I_b|=|I_f|=|I_g|=3$.

$(6)$ $|I_i|+|I_j|+|I_t|\le10$ for each distinct $i,j,t\in\{1,...,s\}$.\\

\textbf{Claim 3}. Let  $\xi_a,\xi_b,\xi_f$ be  distinct elements of $X$, occuring on $\overrightarrow{C}$ in a consecutive order.  If $\xi_a^-\xi_b^+\in E(G)$ then $w\xi_a, w\xi_b \not\in E(G)$ for each $w\in \{\xi_f^+,\xi_f^-\}$.

\textbf{Proof}. If $\xi_f^-\xi_a\in E(G)$ then
$$
\xi_fx_1\xi_b\overleftarrow{C}\xi_a \xi_f^-\overleftarrow{C}\xi_b^+\xi_a^-\overleftarrow{C}\xi_f
$$
is longer than $C$, a contradiction. If $\xi_f^-\xi_b\in E(G)$ then
$$
\xi_fx_1\xi_a\overrightarrow{C}\xi_b \xi_f^-\overleftarrow{C}\xi_b^+\xi_a^-\overleftarrow{C}\xi_f
$$
is longer than $C$, a contradiction.  So, $\xi_{f}^-\xi_a, \xi_{f}^-\xi_b\not\in E(G)$.  By a symmetric argument, $\xi_{f}^+\xi_a, \xi_{f}^+\xi_b\not\in E(G)$.           \  \  \   $\Delta$\\

\textbf{Claim 4}. Let $a,b\in \{1,...,s\}$. If $\xi_a^+w\in E(G)$ for some $w\in V(\xi_b\overrightarrow{C}\xi_a^-)$ then $\xi_b^-y\not\in E(G)$ for each $y\in\{w^+,w^-\}$.

\textbf{Proof}. If $\xi_b^-w^-\in E(G)$ then
$$
\xi_ax_1\xi_b\overrightarrow{C}w^-\xi_b^-\overleftarrow{C}\xi_a^+w\overrightarrow{C}\xi_a
$$
is longer than $C$, a contradiction.  If $\xi_b^-w^+\in E(G)$ then
$$
\xi_ax_1\xi_b\overrightarrow{C}w\xi_a^+\overrightarrow{C}\xi_b^-w^+\overrightarrow{C}\xi_a
$$
is longer than $C$, a contradiction. So, $\xi_b^-y\not\in E(G)$ for each $y\in\{w^+,w^-\}$.      \  \  \   $\Delta$\\

\textbf{Claim 5}. Let $a,b\in \{1,...,s\}$. If $\xi_a^+w\in E(G)$ for some $w\in V(\xi_b^+\overrightarrow{C}\xi_a)$ then $\xi_b^+w^+\not\in E(G)$. If $\xi_a^+w\in E(G)$ for some $w\in V(\xi_a^+\overrightarrow{C}\xi_b)$ then $\xi_b^+w^-\not\in E(G)$. 

\textbf{Proof}. If $\xi_a^+w\in E(G)$ for some $w\in V(\xi_b^+\overrightarrow{C}\xi_a)$ and $\xi_b^+w^+\in E(G)$ then 
$$
\xi_ax_1\xi_b\overleftarrow{C}\xi_a^+w\overleftarrow{C}\xi_b^+w^+\overrightarrow{C}\xi_a
$$ 
is longer than $C$, a contradiction. Hence $\xi_b^+w^+\not\in E(G)$. By a symmetric argument, if  $\xi_a^+w\in E(G)$ for some $w\in V(\xi_a^+\overrightarrow{C}\xi_b)$ then $\xi_b^+w^-\not\in E(G)$.       \  \  \   $\Delta$\\

If $\Upsilon(I_1,...,I_s)=\emptyset$ then $G\backslash \{\xi_1,...,\xi_s\}$ has at least $s+1$ components, contradicting the fact that $\tau>1$. Otherwise   $\Upsilon(I_a,I_b)\not=\emptyset$ for some distinct $a,b\in \{1,...,s\}$.  By Lemma 2, $|I_a|+|I_b|\ge 6$. Since $C$ is extreme, we have  $|I_i|\ge2$ $(i=1,...,s)$. Assume first that   $s\ge\delta+1$. Then   
$$
c= |I_a|+|I_b|+\sum_{i\in\{1,...,s\}\backslash\{a,b\}}|I_i|\ge 6+2(s-2)\ge2\delta+4.
$$
By (1), $c=2\delta+4$ and $|I_a|+|I_b|=6$. By Lemma 2,  $\Upsilon(I_a,I_b)$ consists of a single edge $yz$ with $y\in V(I_a^\ast)$ and $y\in V(I_a^\ast)$. If $|I_a|=|I_b|=3$ then by Lemma 2, $\Upsilon(I_1,...,I_s)=\Upsilon(I_a,I_b)=\{yz\}$ and therefore,  $G\backslash \{\xi_1,...,\xi_s,y\}$ has at least $s+1$ components, contradicting the fact that $\tau>1$. Now let $|I_a|=4$ and $|I_b|=2$. Put $I_a=\xi_aw_1w_2w_3\xi_{a+1}$ and $I_b=\xi_bw_4\xi_{b+1}$. By Claim 1, $y=w_2$ and $z=w_4$. Since $|I_i|=2$ for each $i\in\{1,...,s\}\backslash \{a\}$, we can state that $w_2$ belongs to all edges in  $\Upsilon(I_1,...,I_s)$. Then $G\backslash \{\xi_1,...,\xi_s,w_2\}$ has at least $s+1$ components, again contradicting the fact that $\tau>1$. So, 
$$
s=\delta.                    \eqno{(2)}
$$
Recalling that $\Upsilon(I_a,I_b)\not=\emptyset$, we can choose $L\in \Upsilon(I_a,I_b)$.            If $|L|\ge3$ then by Lemma 2, 
$$
|I_a|+|I_b|\ge2\overline{p}+2|L|+4\ge10,
$$
contradicting Claim 2(1). Otherwise $|L|\le2$. \\

\textbf{Claim 6}. $\Upsilon(I_1,...,I_s)$ consists of pairwise edge disjoint paths of length at most two.

\textbf{Proof}. Assume the contrary, that is $P_1,P_2\in \Upsilon(I_1,...,I_s)$ and $P_1=y_1y_2y_3$, $P_2=y_1y_2y_4$. If $y_1,y_2,y_4$ belong to different elementary segments $I_a,I_b,I_f$ then  by Lemma 2,
$$
|I_a|+|I_b|\ge8, \  \ |I_a|+|I_f|\ge8, \  \   |I_b|+|I_f|\ge8. 
$$
This implies  $|I_a|+|I_b|+|I_f|\ge12$, contradicting Claim 2(6). Now let $y_1\in V(I_a^\ast)$ and $y_3,y_4\in V(I_b^\ast)$. Assume w.l.o.g. that $y_3\in V(\xi_b^+\overrightarrow{C}y_4^-)$. Put
$$
|\xi_a\overrightarrow{C}y_1|=d_1,  \  \  |y_1\overrightarrow{C}\xi_{a+1}|=d_2,
$$
$$
|\xi_b\overrightarrow{C}y_3|=d_3,   \  \ |y_3\overrightarrow{C}y_4|=d_4,  \  \  |y_4\overrightarrow{C}\xi_{b+1}|=d_5.
$$
Since $C$ is extreme, we have
$$
|C|\ge|\xi_ax_1\xi_b\overleftarrow{C}y_1y_2y_3\overrightarrow{C}\xi_a|\ge|C|-d_1-d_3+4,
$$
$$
|C|\ge|\xi_a\overrightarrow{C}y_1y_2y_4\overleftarrow{C}\xi_{a+1}x_1\xi_{b+1}|\ge|C|-d_2-d_5+4,
$$
implying that $d_1+d_3\ge4$ and $d_2+d_5\ge4$. Observing also that $d_4\ge|y_3y_2y_4|=2$, we have $|I_a|+|I_b|\ge\sum_{i=1}^5d_i\ge10$, contradicting Claim 2(1).     \  \  \      $\Delta$\\

By Claim 2(4), $|i:|I_i|\ge3|\le 4$. Further, if $|i:|I_i|\ge3|=0$ then by Lemma 2, $\Upsilon(I_1,...,I_s)=\emptyset$, implying that $\tau<1$, a contradiction. So, 
$$
1\le|i:|I_i|\ge3|\le 4.
$$

\textbf{Case 1.1}. $|i:|I_i|\ge3|=4$.

Let $|I_i|\ge3$ $(i=a,b,f,g)$ and $|I_i|=2$ for each $i\in\{1,...,s\}\backslash \{a,b,f,g\}$. By Claim 2(5), $|I_i|=3$ $(i=a,b,f,g)$. Assume w.l.o.g. that $\xi_a,\xi_b,\xi_f,\xi_g$ occur on $C$ in a consecutive order. By Lemma 2, $\Upsilon(I_1,...,I_s)\subseteq E(G)$.\\

 \textbf{Claim 7}. If $\Upsilon(I_a,I_f)\not=\emptyset$ then $\Upsilon(I_b,I_g)=\emptyset$.

\textbf{Proof}. Assume the contrary, that is $\Upsilon(I_b,I_g)\not=\emptyset$. Let $y_1y_2\in \Upsilon(I_a,I_f)$, where $y_1\in V(I_a^\ast)$ and $y_2\in V(I_f^\ast)$. Assume w.l.o.g. $y_1=w_2$. By Claim 1, $y_2=w_5$. Analogously, there is an intermediate edge $y_3y_4$ between $I_b^\ast$ and $I_g^\ast$. Assume w.l.o.g. that $y_3=w_3$. By Claim 1, $y_4=w_8$. Further, we have $N(w_7)\cap \{\xi_{g+1}, \xi_{a+1}, \xi_b,\xi_f\}=\emptyset$ (by Claim 3), $N(w_7)\cap \{w_4,w_6\}=\emptyset$ (by Claim 4), $N(w_7)\cap \{w_1,w_3,w_5\}=\emptyset$ (by Claim 1).  If $N(w_7)\subseteq V(C)$ then
$$
N(w_7)\subseteq \{\xi_1,...,\xi_s,w_8,w_2\}\backslash \{\xi_{a+1},\xi_b,\xi_f,\xi_{g+1}\},
$$  
implying that $|N(w_7)|\le s-1=\delta-1$, a contradiction. Now let $N(w_7)\not\subseteq V(C)$, that is $x_2w_7\in E(G)$ for some $x_2\in V(G\backslash C)$. Clearly, $x_2\not=x_1$, $N(x_2)\subseteq V(C)$ and $x_2x_8, x_2\xi_g\not\in E(G)$. Observing also that $x_2w_2\not\in E(G)$  (by Lemma 2), we have
$$
N(w_7)\subseteq \{\xi_1,...,\xi_s,w_7\}\backslash \{\xi_{a+1},\xi_b,\xi_f,\xi_g\},
$$
a contradiction.     \  \  \  $\Delta$              \\

\textbf{Claim 8}. $\Upsilon(I_a,I_f)=\Upsilon(I_b,I_g)=\emptyset$.

\textbf{Proof}. Assume the contrary, that is $\Upsilon(I_a,I_f)\not=\emptyset$. As in proof of Claim 7, assume w.l.o.g. that $w_2w_5\in E(G)$. Then we have $N(w_7)\cap \{w_1,w_3,w_5\}=\emptyset$ (by Claim 1), $N(w_7)\cap \{\xi_{a+1},\xi_f\}=\emptyset$ (by Claim 3),   $w_7w_6\not\in E(G)$ (by Claim 4) and  $w_7w_4\not\in E(G)$ (by Claim 7). If $N(w_7)\not\subseteq V(C)$ then  $x_2w_7\in E(G)$ for some $x_2\in V(G\backslash C)$. Clearly, $x_2\not=x_1$, $N(x_2)\subseteq V(C)$ and $x_2x_8, x_2\xi_g\not\in E(G)$. Observing also that $x_2w_2\not\in E(G)$ (by Lemma 2), we have
$$
N(x_2)\subseteq \{\xi_1,...,\xi_s,w_7\}\backslash \{\xi_{a+1},\xi_f,\xi_g\},
$$
contradicting the fact that $|N(x_2)|\ge\delta=s$. Now let $N(w_7)\subseteq V(C)$. Then
$$
N(w_7)\subseteq\{\xi_1,...,\xi_s,w_8,w_2\}\backslash \{\xi_{a+1},\xi_f\},
$$
implying that $w_7w_2\in E(G)$. But then $N(w_1)\cap \{w_3,w_5,w_7\}=\emptyset$ (by Claim 1),  $N(w_1)\cap \{\xi_{a+1},\xi_f,\xi_g\}=\emptyset$ (by Claim 3), $w_1w_6\not\in E(G)$ (by Claim 4) and  $w_1,w_8\not\in E(G)$ (by Lemma 2). If $N(w_1)\not\subseteq V(C)$ then $x_3w_1\in E(G)$ for some $x_3\in V(G\backslash C)$. Clearly, $x_3\not=x_1$, $N(x_3)\subseteq V(C)$ and $x_3w_2, x_3\xi_a\not\in E(G)$. Observing also that $x_3w_4\not\in E(G)$ (by Lemma 2), we have
$$
N(x_3)\subseteq \{\xi_1,...,\xi_s,w_1\}\backslash \{\xi_{a+1},\xi_f,\xi_g,\xi_a\},
$$
contradicting the fact that $|N(x_3)|\ge\delta=s$. Now let $N(w_1)\subseteq V(C)$. Then
$$
N(w_1)\subseteq \{\xi_1,...,\xi_s,w_2,w_4\}\backslash \{\xi_{a+1},\xi_f,\xi_g\},
$$
implying that $|N(w_1)|\le s-1=\delta$, a contradiction.  \  \  \  $\Delta$  \\

\textbf{Claim 9}. $\Upsilon(I_a,I_b,I_f,I_g)=\emptyset$.

\textbf{Proof}. Assume the contrary, that is 
$\Upsilon(I_a,I_b,I_f,I_g)\not=\emptyset$. 
By Claim 8, either $\Upsilon(I_a,I_b)\not=\emptyset$ 
or $\Upsilon(I_b,I_f)\not=\emptyset$ or 
$\Upsilon(I_f,I_g)\not=\emptyset$ or $\Upsilon(I_g,I_a)\not=\emptyset$.
 Assume w.l.o.g. that $\Upsilon(I_a,I_g)\not=\emptyset$.
 By Claim 1, either $w_2w_7\in E(G)$ or $w_1w_8\in E(G)$.\\

\textbf{Case a}. $w_2w_7\in E(G)$.

We have $w_8w_2, w_8w_6\not\in E(G)$ (by Claim 1), 
  $w_8\xi_{a+1}, w_8\xi_g\not\in E(G)$ (by Claim 3), 
$w_8w_3, w_8w_4\not\in E(G)$ (by Claim 8) and 
 $w_8w_1\not\in E(G)$ (by Lemma 2). If   $N(w_8)\not\subseteq V(C)$ then $x_2w_8\in E(G)$ for some $x_2\in V(G\backslash C)$. Clearly, $x_2\not=x_1$, $N(x_2)\subseteq V(C)$ and $x_2w_7, x_2\xi_{g+1}\not\in E(G)$. Observing also that $x_2w_5\not\in E(G)$ (by Lemma 2), we have
$$
N(x_2)\subseteq \{\xi_1,...,\xi_s,w_8\}\backslash \{\xi_{a+1},\xi_g,\xi_{g+1}\},
$$
contradicting the fact that $|N(x_2)|\ge\delta=s$. Now let $N(w_8)\subseteq V(C)$. Then
$$
N(w_8)\subseteq\{\xi_1,...,\xi_s,w_5,w_7\}\backslash \{\xi_{a+1},\xi_g\},
$$
implying that $w_8w_5\in E(G)$. Then we have $w_1w_3, w_1w_7\not\in E(G)$ (by Claim 1), 
$w_1\xi_{a+1},w_1\xi_g\not\in E(G)$ (by Claim 3), $w_1w_5, w_1w_6\not\in E(G)$ 
(by Claim 8) and  $w_1w_8\not\in E(G)$ (by Lemma 2). If   $N(w_1)\not\subseteq V(C)$ then $x_3w_1\in E(G)$ for some $x_3\in V(G\backslash C)$. Clearly, $x_3\not=x_1$, $N(x_3)\subseteq V(C)$ and $x_3w_2, x_3\xi_a\not\in E(G)$. Observing also that $x_3w_4\not\in E(G)$ (by Lemma 2), we have
$$
N(x_3)\subseteq \{\xi_1,...,\xi_s,w_1\}\backslash \{\xi_{a+1},\xi_g,\xi_a\},
$$
contradicting the fact that $|N(x_2)|\ge\delta=s$. Now let $N(w_1)\subseteq V(C)$. Then 
$$
N(w_1)\subseteq \{\xi_1,...,\xi_s,w_2,w_4\}\backslash \{\xi_{a+1},\xi_g\},
$$
implying that $w_1w_4\in E(G)$. By a symmetric argument, $w_6w_3\in E(G)$. But then
$$
\xi_aw_1w_4w_3w_6\overrightarrow{C}w_7w_2\overrightarrow{C}\xi_bx_1\xi_{b+1}\overrightarrow{C}w_5w_8\overrightarrow{C}\xi_a
$$
is longer than $C$, a contradiction.\\

\textbf{Case b}. $w_1w_8\in E(G)$.

If $\Upsilon(I_a,I_b,I_f,I_g)=\{w_1w_8\}$ then 
$G\backslash \{\xi_1,...,\xi_s,w_1\}$ has at least $s+1$ components, contradicting the fact that $\tau>1$.
Let $\Upsilon(I_a,I_b,I_f,I_g)\not=\{w_1w_8\}$. 
If either $w_1w_4\in E(G)$ or $w_3w_6\in E(G)$ 
or $w_5w_8\in E(G)$ then we can argue as in Case a. 
Otherwise, by Claim 8, either $w_2w_3\in E(G)$ 
or $w_4w_5\in E(G)$ or $w_6w_7\in E(G)$. 
Observing that $w_2w_3, w_6w_7\not\in E(G)$ (by Claim 4), 
we have $w_4w_5\in E(G)$. Then we have 
$w_2w_4\not\in E(G)$ (by Claim 1), 
$w_2\xi_a, w_2\xi_{b+1}, w_2\xi_f, w_2\xi_{g+1}\not\in E(G)$
 (by Claim 3), $w_2w_3\not\in E(G)$ (by Claim 4), 
$w_2w_5, w_2w_6\not\in E(G)$ (by Claim 8) and  
$w_2w_7, w_2w_8\not\in E(G)$ (by Lemma 2). If $N(w_2)\not\subseteq V(C)$ then $x_2w_2\in E(G)$ for some $x_2\in V(G\backslash C)$. Clearly, $x_2\not=x_1$, $N(x_2)\subseteq V(C)$ and $x_2w_1, x_2\xi_{a+1}\not\in E(G)$. Then
$$
N(x_2)\subseteq \{\xi_1,...,\xi_s,w_2\}\backslash \{\xi_a,\xi_{a+1},\xi_{b+1},\xi_f,\xi_{g+1}\},
$$
contradicting the fact that $|N(x_2)|\ge\delta=s$. Now let $N(w_2)\subseteq V(C)$. Then
$$
N(w_2)\subseteq \{\xi_1,...,\xi_s,w_1\}\backslash \{\xi_a, \xi_{b+1},\xi_f, \xi_{g+1}\},
$$ 
a contradiction.\\

\textbf{Case 1.2}. $|i:|I_i|\ge3|=3$.

Let $|I_i|\ge3$ $(i=a,b,f)$ and $|I_i|=2$ for each $i\in\{1,...,s\}\backslash \{a,b,f\}$. By Claim 2(6), $9\le|I_a|+|I_b|+|I_f|\le10$.\\

 \textbf{Case 1.2.1}.  $|I_a|+|I_b|+|I_f|=10$.

Assume w.l.o.g. that $|I_a|=|I_b|=3$ and $|I_f|=4$. Put 
$$
I_a=\xi_aw_1w_2\xi_{a+1},  \  \  I_b=\xi_bw_3w_4\xi_{b+1},  \  \  I_f=\xi_fw_5w_6w_7\xi_{f+1}  
$$
By Lemma 2, $\Upsilon(I_1,...,I_s)\subseteq E(G)$.\\

\textbf{Case 1.2.1.1}. $\Upsilon(I_a,I_b)\not=\emptyset$.

By Claim 1, either $w_2w_3\in E(G)$ or $w_1w_4\in E(G)$.\\

\textbf{Case 1.2.1.1.1}.  $w_2w_3\in E(G)$.

If $\Upsilon(I_1,...,I_s)=\{w_2w_3\}$ then $G\backslash \{\xi_1,...,\xi_s,w_2\}$ has at lest $s+1$ components, contradicting the fact that $\tau>1$. Let $\Upsilon(I_1,...,I_s)\not=\{w_2w_3\}$. Further, if $\Upsilon(I_f,I_a)=\Upsilon(I_f,I_b)=\emptyset$ then by Claim 1, $w_6$ belongs to every edge in $\Upsilon(I_1,...,I_s)\backslash\{w_2w_3\}$ connecting $I_f$ with some segment of length 2. By Claim 5, $w_5w_7\not\in E(G)$ and hence $G\backslash \{\xi_1,...,\xi_s,w_6,w_3\}$ has at least $s+2$ components, a contradiction. Now let either $\Upsilon(I_f,I_a)\not=\emptyset$ or $\Upsilon(I_f,I_b)\not=\emptyset$, say $\Upsilon(I_f,I_b)\not=\emptyset$. By Claim 1, either $w_4w_6\in E(G)$ or $w_3w_6\in E(G)$ or $w_3w_7\in E(G)$.\\

\textbf{Case 1.2.1.1.1.1}.  $w_4w_6\in E(G)$.

We have $w_5w_1, w_4w_7, w_7w_2, w_4w_2\not\in E(G)$ (by Claim 1), $w_5w_4, w_7w_1\not\in E(G)$ (by Claim 4), $w_5w_2, w_5w_7\not\in E(G)$ (by Claim 5),  $ w_4w_1\not\in E(G)$ (by Lemma 2).  Then $G\backslash \{\xi_1,...,\xi_s,w_3,w_6\}$ has at least $s+2$ components, a contradiction.\\

\textbf{Case 1.2.1.1.1.2}.  $w_3w_7\in E(G)$.

We have $w_4w_2, w_4w_7\not\in E(G)$ (by Claim 1), $w_4\xi_{f+1}, w_4\xi_{a+1}, w_4\xi_b\not\in E(G)$ (Claim 3) and  $w_4w_1,w_4w_5,w_4w_6\not\in E(G)$ (by Lemma 2).  If $N(w_4)\not\subseteq V(C)$ then  $x_2w_4\in E(G)$ for some $x_2\in V(G\backslash C)$. Clearly, $x_2\not=x_1$, $N(x_2)\subseteq V(C)$ and $x_2w_3, x_2\xi_{b+1}\not\in E(G)$. Then
$$
N(x_2)\subseteq \{\xi_1,...,\xi_s,w_4\}\backslash \{\xi_{a+1},\xi_b,\xi_{b+1},\xi_{f+1}\},
$$
a contradiction. Now let  $N(w_4)\subseteq V(C)$. Then 
$$
N(w_4)\subseteq \{\xi_1,...,\xi_s,w_3\}\backslash \{\xi_{f+1},\xi_{a+1},\xi_b\},
$$
a contradiction.\\

\textbf{Case 1.2.1.1.1.3}.  $w_3w_6\in E(G)$.

We have $w_1w_3,w_1w_5,w_3w_5,w_4w_7\not\in E(G)$ (by Claim 1), $w_1w_7\not\in E(G)$ (Claim 4), $w_3w_7\not\in E(G)$ (otherwise we can argue as in Case 1.2.1.1.1.2),  $w_5w_7\not\in E(G)$ (Claim 5) and  $w_1w_4,w_4w_5\not\in E(G)$ (by Lemma 2). So, $G\backslash \{\xi_1,...,\xi_s,w_2,w_6\}$ has at least $s+2$ components, a contradiction.\\

\textbf{Case 1.2.1.1.2}.  $w_1w_4\in E(G)$.

We have $w_3w_1,w_3w_5\not\in E(G)$ (by Claim 1), $w_3w_2\not\in E(G)$ (by Lemma 2) and  $w_3\xi_a,w_3\xi_{b+1}\not \in E(G)$ (by Claim 3). Assume first that $N(w_3)\not\subseteq V(C)$, that is   $x_2w_3\in E(G)$ for some $x_2\in V(G\backslash C)$. Clearly, $x_2\not=x_1$, $N(x_2)\subseteq V(C)$ and $x_2w_4, x_2\xi_b\not\in E(G)$. Observing also that $x_2w_6,x_2w_7\not\in E(G)$ (by Lemma 2), we have
$$
N(x_2)\subseteq \{\xi_1,...,\xi_s,w_3\}\backslash \{\xi_a,\xi_{b+1}\},
$$
a contradiction. Now let $N(w_3)\subseteq V(C)$. If $N(w_3)\cap \{w_6,w_7\}=\emptyset$ then 
$$
N(w_3)\subseteq \{\xi_1,...,\xi_s,w_4\}\backslash \{\xi_a,\xi_{b+1}\},
$$
a contradiction. Hence $N(w_3)\cap \{w_6,w_7\}\not=\emptyset$. By a symmetric argument,  $N(w_2)\cap \{w_5,w_6\}\not=\emptyset$. We have three main subcases, namely either $w_2w_5$, $w_3w_7\in E(G)$ or $w_2w_6,w_3w_6\in E(G)$ or $w_2w_6,w_3w_7\in E(G)$.\\

\textbf{Case 1.2.1.1.2.1}.  $w_2w_5,w_3w_7\in E(G)$.

If $w_2w_6\in E(G)$ then
$$
\xi_aw_1w_4\overleftarrow{C}\xi_{a+1}x_1\xi_{b+1}\overrightarrow{C}w_5w_2w_6\overrightarrow{C}\xi_a
$$
is longer than $C$, a contradiction. Hence, $w_2w_6\not\in E(G)$. By a symmetric argument, $w_3w_6\not\in E(G)$. We have $w_6w_1, w_6w_4\not\in E(G)$ (by Lemma 2). If $w_6\xi_{a+1}\in E(G)$ then 
$$
\xi_ax_1\xi_{b+1}\overrightarrow{C}w_5w_2w_1w_4\overleftarrow{C}\xi_{a+1}w_6\overrightarrow{C}\xi_a
$$
is longer than $C$, a contradiction. Let $w_6\xi_{a+1}\not\in E(G)$. By a symmetric argument, $w_6\xi_b\not\in E(G)$. Further, if $w_6\xi_a\in E(G)$ then
$$
\xi_aw_6w_5w_2w_1w_4\overrightarrow{C}\xi_fx_1\xi_{a+1}\overrightarrow{C}w_3w_7\overrightarrow{C}\xi_a
$$
is longer than $C$, a contradiction. Hence, $w_6\xi_a\not\in E(G)$. By a symmetric argument, $w_6\xi_{b+1}\not\in E(G)$. Assume that $N(w_6)\not\subseteq V(C)$, that is $x_2w_6\in E(G)$ for some $x_2\in V(G\backslash C)$. Clearly, $x_2\not=x_1$, $N(x_2)\subseteq V(C)$ and $x_2w_5, x_2w_7\not\in E(G)$. Then we have
$$
N(x_2)\subseteq \{\xi_1,...,\xi_s,w_6\}\backslash \{\xi_a,\xi_{a+1},\xi_b,\xi_{b+1}\},
$$
a contradiction. Now let $N(w_6)\subseteq V(C)$. Then
$$
N(w_6)\subseteq \{\xi_1,...,\xi_s,w_5,w_7\}\backslash \{\xi_a,\xi_{a+1},\xi_b,\xi_{b+1}\},
$$
a contradiction.\\

\textbf{Case 1.2.1.1.2.2}.  $w_2w_6,w_3w_6\in E(G)$.

By Claim 5, $w_5w_7\not\in E(G)$. If $w_2w_3\in E(G)$ then we can argue as in Case 1.2.1.1.1. Let $w_2w_3\not\in E(G)$. Next, if $w_2w_5\in E(G)$ then we can argue as in Case 1.2.1.1.2.1. Let $w_2w_5\not\in E(G)$. We have also $w_2w_7,w_3w_5\not\in E(G)$ (by Claim 1), $w_3w_7\not\in E(G)$ (as in Case 1.2.1.1.2.1). Then $G\backslash \{\xi_1,...,\xi_s,w_4,w_6\}$ has at least $s+2$ components, contradicting the fact that $\tau>1$.\\

\textbf{Case 1.2.1.1.2.3}.  $w_2w_6,w_3w_7\in E(G)$.

We have $w_5w_1, w_5w_3\not\in E(G)$ (by Claim 1), $w_5w_2\not\in E(G)$ (as in Case 1.2.1.1.2.1), $w_5\xi_b,w_5\xi_{f+1}\not\in E(G)$ (by Claim 3), $w_5w_4\not\in E(G)$ (by Lemma 2), $w_5w_7\not\in E(G)$ (as in Case 1.2.1.1.2.2). Assume that $N(w_5)\not\subseteq V(C)$, that is  $x_2w_5\in E(G)$ for some $x_2\in V(G\backslash C)$. Clearly, $x_2\not=x_1$, $N(x_2)\subseteq V(C)$ and $x_2w_6, x_2\xi_f\not\in E(G)$. Then we have
$$
N(x_2)\subseteq \{\xi_1,...,\xi_s,w_5,w_7\}\backslash \{\xi_{f+1},\xi_b,\xi_f\},
$$
a contradiction. Now let $N(w_5)\subseteq V(C)$. Then  
$$
N(w_5)\subseteq \{\xi_1,...,\xi_s,w_6\}\backslash \{\xi_{f+1},\xi_b\},
$$
a contradiction.\\

\textbf{Case 1.2.1.2}. $\Upsilon(I_a,I_b)=\emptyset$.

It follows that $w_6$ belongs to all edges in $\Upsilon(I_1,...,I_s)$, implying that $\tau\le1$, a contradiction.\\

\textbf{Case 1.2.2}. $|I_a|+|I_b|+|I_f|=9$.

It follows that $|I_a|=|I_b|=|I_f|=3$. Put 
$$
I_a=\xi_aw_1w_2\xi_{a+1},  \  \  I_b=\xi_bw_3w_4\xi_{b+1},  \  \  I_f=\xi_fw_5w_6\xi_{f+1}.  
$$
By Lemma 2, $\Upsilon(I_1,...,I_s)=\Upsilon(I_a,I_b,I_f)$. Assume w.l.o.g. that $\Upsilon(I_a,I_b)\not=\emptyset$ and $w_2w_3\in E(G)$. We have $w_1w_6,w_4w_5\not\in E(G)$ (Claim 4). If $\Upsilon(I_1,...,I_s)=\Upsilon(I_a,I_b)$ then clearly $\tau\le1$, a contradiction. Otherwise, assume w.l.o.g. that  $w_2w_5\in E(G)$. Hence, $w_1\xi_{a+1},w_1\xi_b,w_1\xi_f\not\in E(G)$ (by Claim 3) and  $w_1w_4$, $w_1w_5$, $w_1w_6\not\in E(G)$. Assume that $N(w_1)\not\subseteq V(C)$, that is  $x_2w_1\in E(G)$ for some $x_2\in V(G\backslash C)$. Clearly, $x_2\not=x_1$, $N(x_2)\subseteq V(C)$ and $x_2w_2, x_2\xi_a\not\in E(G)$. Then we have
$$
N(x_2)\subseteq \{\xi_1,...,\xi_s,w_4\}\backslash \{\xi_a,\xi_{a+1},\xi_b,\xi_f\},
$$
a contradiction. Now let $N(w_1)\subseteq V(C)$. Then 
$$
N(w_1)\subseteq \{\xi_1,...,\xi_s,w_2\}\backslash \{\xi_{a+1},\xi_b,\xi_f\},
$$ 
a contradiction.\\

\textbf{Case 1.3}. $|i:|I_i|\ge3|=2$.

Let $|I_i|\ge3$ $(i=a,b)$ and $|I_i|=2$ for each $i\in\{1,...,s\}\backslash \{a,b\}$. By Claim 2(2), $6\le|I_a|+|I_b|\le8$.\\

\textbf{Case 1.3.1}. $|I_a|+|I_b|=8$.

\textbf{Case 1.3.1.1}. $|I_a|=3$, $|I_b|=5$.

Put $I_a=\xi_aw_1w_2\xi_{a+1}$ and $I_b=\xi_bw_3w_4w_5w_6\xi_{b+1}$.\\

\textbf{Case 1.3.1.1.1}. $\Upsilon(I_1,...,I_s)=\Upsilon(I_a,I_b)$.

If there is a vertex belonging to all edges in $\Upsilon(I_a,I_b)$  then clearly $\tau\le1$, a contradiction. Otherwise $w_1y_1, w_2y_2\in E(G)$ for some distinct $y_1,y_2\in\{w_3,w_4,w_5,w_6\}$. By Claim 4, $y_1,y_2$ are not consequent vertices on $C$ and $\{y_1,y_2\}\not=\{w_3,w_6\}$. Then we can assume w.l.o.g. that $w_2w_3,w_1w_5\in E(G)$. If $w_4w_6\in E(G)$ then 
$$
\xi_ax_1\xi_{a+1}\overrightarrow{C}w_3w_2w_1w_5w_4w_6\overrightarrow{C}\xi_a
$$
is longer than $C$, a contradiction. Let $w_4w_6\not\in E(G)$. Then $G\backslash \{\xi_1,...,\xi_s,w_3,w_5\}$ has at least $s+2$ components, a contradiction.\\

\textbf{Case 1.3.1.1.2}. $\Upsilon(I_1,...,I_s)\not=\Upsilon(I_a,I_b)$.

Choose a segment $I_f=\xi_fw_7\xi_{f+1}$ such that $w_7y\in E(G)$ for some $y\in V(I_a^\ast)\cup V(I_b^\ast)$. By Lemma 2, $y\in V(I_b^\ast)$. Assume w.l.o.g. that $\xi_a,\xi_b,\xi_f$ occur on $C$ in this order. By Claim 1, either $y=w_4$ or $y=w_5$.\\

\textbf{Case 1.3.1.1.2.1}. $y=w_4$.

Let $g\in\{1,...,s\}\backslash\{a,b\}$. Clearly $|I_g|=2$. Put $I_g=\xi_gw_8\xi_{g+1}$. By Claim 4, $w_7w_5,w_8w_5\not\in E(G)$, i.e. $N(w_8)\cap \{w_3,w_4,w_5,w_6\}\subseteq\{w_4\}$. Further, we have $N(w_2)\cap\{w_3,w_5\}=\emptyset$ (by Claim 4) and $w_1w_5\not\in E(G)$ (by claim 5). If $w_3w_5\in E(G)$ then
$$
\xi_a\overrightarrow{C}\xi_bx_1\xi_f\overleftarrow{C}w_5w_3w_4w_7\overrightarrow{C}\xi_a
$$
is longer than $C$. Let $w_3w_5\not\in E(G)$. So, $G\backslash\{\xi_1,...,\xi_s,w_4,w_6\}$ has at least $s+2$ components, a contradiction.\\

\textbf{Case 1.3.1.1.2.2}. $y=w_5$.

By Claim 5,  $w_4w_6, w_3w_6\not\in E(G)$. In addition, we have $w_1w_3,w_2w_6\not\in E(G)$ (by Claim 1), $w_2w_4\not\in E(G)$ (by Claim 4) and $w_1w_6\not\in E(G)$ (by Claim 5). If $w_1w_4\not\in E(G)$ then $G\backslash \{\xi_1,...,\xi_s,w_3,w_5\}$ has at least $s+2$ components, a contradiction. Now let $w_1w_4\in E(G)$, implying that $w_1w_3\not\in E(G)$ (by Claim 1) and $w_2w_3\not\in E(G)$ (by Claim 4).  But then $G\backslash \{\xi_1,...,\xi_s,w_4,w_5\}$ has at least $s+2$ components, again a contradiction.\\

\textbf{Case 1.3.1.2}. $|I_a|=|I_b|=4$.

Put $I_a=\xi_aw_1w_2w_3\xi_{a+1}$ and $I_b=\xi_bw_4w_5w_6\xi_{b+1}$.\\

\textbf{Case 1.3.1.2.1}. $\Upsilon(I_1,...,I_s)\not=\Upsilon(I_a,I_b)$.

Assume w.l.o.g. that $\Upsilon(I_a,I_f)\not=\emptyset$ for some $f\in\{1,...,s\}\backslash\{a,b\}$, and $\xi_a,\xi_b,\xi_f$ occur on $C$ in this order. Put $I_f=\xi_fw_7\xi_{f+1}$. By Claim 1, $w_7w_2\in E(G)$ and by Claim 5,  $w_1w_3\not\in E(G)$.  Observing also that $w_1w_4,w_3w_6\not\in E(G)$ (by Claim 1), $w_3w_4\not\in E(G)$ (by Claim 4) and $w_1w_6\not\in E(G)$ (by Claim 5), we conclude that $G\backslash \{\xi_1,...,\xi_s,w_2,w_5\}$ has at least $s+2$ components, a contradiction.\\

\textbf{Case 1.3.1.2.2}. $\Upsilon(I_1,...,I_s)=\Upsilon(I_a,I_b)$.

If there is a vertex belonging to all edges in $\Upsilon(I_1,...,I_s)$  then $\tau\le1$, a contradiction. Otherwise, by Claim 4, either $w_3w_4,w_2w_5\in E(G)$ or $w_3w_5,w_2w_6\in E(G)$ or $w_3w_4,w_1w_6\in E(G)$ or $w_3w_5,w_2w_4\in E(G)$.\\

\textbf{Case 1.3.1.2.2.1}. $w_3w_4,w_2w_5\in E(G)$.

If $w_1w_3\in E(G)$ then
$$
\xi_ax_1\xi_{a+1}\overrightarrow{C}w_4w_3w_1w_2w_5\overrightarrow{C}\xi_a
$$ 
is longer than $C$, a contradiction. Let $w_1w_3\not\in E(G)$. Next, if $w_1w_6\in E(G)$ then
$$
\xi_ax_1\xi_b\overleftarrow{C}w_3w_4w_5w_2w_1w_6\overrightarrow{C}\xi_a
$$
is longer than $C$, a contradiction. Let $w_1w_6\not\in E(G)$. Observe that $w_1w_4\not\in E(G)$ (by Claim 1), $w_1\xi_{a+1},w_1\xi_b\not\in E(G)$ (by Claim 3) and   $w_1w_5\not\in E(G)$ (Claim 4). 
 Moreover, if $f\in\{1,...,s\}\backslash\{a,b\}$ and $I_f=\xi_fw_7\xi_{f+1}$ 
then by Claim 1, $w_1w_7\not\in E(G)$.
 So, $N(w_1)\subseteq\{\xi_1,...,\xi_s,w_2\}\backslash \{\xi_{a+1},\xi_b\}$. 
If $\xi_{a+1}\not=\xi_b$ then $|N(w_1)|\le s-1=\delta-1$, 
a contradiction. Let $\xi_{a+1}=\xi_b$. Since $s\ge3$, 
we have $\xi_{b+1}\not=\xi_a$. Further, we have 
$w_7w_1,w_7w_3\not\in E(G)$ (by Claim 1), $w_7\xi_{a+1}\not\in E(G)$ (by Claim 3) and $w_7w_2\not\in E(G)$  (by Claim 4).
 Hence, if $N(w_7)\subseteq V(C)$ then $N(w_7)\subseteq\{\xi_1,...,\xi_s\}\backslash\{\xi_{a+1}\}$, a contradiction. Analogous arguments can be used when $N(w_7)\not\subseteq V(C)$.   \\

\textbf{Case 1.3.1.2.2.2}. $w_3w_5,w_2w_6\in E(G)$.

By Claim 5, If $w_1w_3\not\in E(G)$.     By a symmetric argument, $w_4w_6\not\in E(G)$. We have also $w_1w_4,w_3w_6\not\in E(G)$ (by Claim 1) and  $w_1w_6,w_3w_4\not\in E(G)$ (by Claim 4).  So, $G\backslash \{\xi_1,...,\xi_s,w_2,w_5\}$ has at least $s+3$ components, a contradiction.\\

\textbf{Case 1.3.1.2.2.3}. $w_3w_4,w_1w_6\in E(G)$.

If $\xi_a=\xi_{b+1}$ and $\xi_{a+1}=\xi_b$ then clearly $s=2$, a contradiction. Assume w.l.o.g. that $\xi_a\not=\xi_{b+1}$. Choose $f\in\{1,...,s\}\backslash \{a,b\}$ such that  $\xi_a,\xi_b,\xi_f$ occur on $C$ in this order. Clearly, $|I_f|=2$. Put $I_f=\xi_fw_7\xi_{f+1}$. Then $w_7w_1,w_7w_3,w_7w_4,w_7w_6\not\in E(G)$ (Claim 1), $w_7\xi_{a+1},w_7\xi_b\not\in E(G)$ (by Claim 3) and $w_7w_2,w_7w_5\not\in E(G)$ (by Claim 5). If $N(w_7)\subseteq V(C)$ then $N(w_7)\subseteq\{\xi_1,...,\xi_s\}\backslash \{\xi_{a+1},\xi_b\}$, a contradiction. Similar arguments can be used when $N(w_7)\not\subseteq V(C)$.\\

\textbf{Case 1.3.1.2.2.4}. $w_3w_5,w_2w_4\in E(G)$.

By Claim 5,  $w_1w_3\not\in E(G)$.    By a similar argument, $w_4w_6\not\in E(G)$. Observing also that $w_1w_4,w_3w_6\not\in E(G)$ (by Claim 1) and $w_1w_6\not\in E(G)$ (by Claim 4), we conclude that if $w_3w_4\not\in E(G)$ then $G\backslash \{\xi_1,...,\xi_s,w_2,w_5\}$ has at least $s+3$ components, a contradiction. Now let $w_3w_4\in E(G)$. Then $w_1\xi_{a+1},w_1\xi_b\not\in E(G)$ (by Claim 3), $w_1w_4\not\in E(G)$ (by Claim 1), $w_1w_5,w_1w_6\not\in E(G)$ (by Claim 4). Assume that $N(w_1)\not\subseteq V(C)$, that is  $x_2w_1\in E(G)$ for some $x_2\in V(G\backslash C)$. Clearly, $x_2\not=x_1$, $N(x_2)\subseteq V(C)$ and $x_2w_2, x_2\xi_a\not\in E(G)$. Then we have
$$
N(x_2)\subseteq \{\xi_1,...,\xi_s,w_1\}\backslash \{\xi_a,\xi_{a+1},\xi_b\},
$$
a contradiction. Now let $N(w_1)\subseteq V(C)$. Then  
$$
N(w_1)\subseteq\{\xi_1,...,\xi_s,w_2\}\backslash \{\xi_{a+1},\xi_b\},
$$
implying that $\xi_{a+1}=\xi_b$. Since $s\ge3$, we have $\xi_a\not=\xi_{b+1}$. Put $I_{a-1}=\xi_{a-1}w_7\xi_a$. We have $w_7w_3,w_7w_4,w_7w_6\not\in E(G)$ (by Claim 1), $w_7\xi_{a+1}\not\in E(G)$ (by Claim 3) and  $w_7w_2,w_7w_5\not\in E(G)$ (by Claim 4).  If $N(w_7)\subseteq V(C)$ then $N(w_7)\subseteq\{\xi_1,...,\xi_s\}\backslash \{\xi_{a+1}\}$, a contradiction. Analogous arguments can be used when $N(w_7)\not\subseteq V(C)$.\\

\textbf{Case 1.3.2}. $|I_a|+|I_b|=7$.

Assume w.l.o.g. that $|I_a|=3$ and $|I_b|=4$. Put $I_a=\xi_aw_1w_2\xi_{a+1}$ and $I_b=\xi_bw_3w_4w_5\xi_{b+1}$. If $\Upsilon(I_a,I_b)=\emptyset$ then $w_4$ belongs to all edges in $\Upsilon(I_1,...,I_s)$, implying that $\tau\le1$, a contradiction. Let $\Upsilon(I_a,I_b)\not=\emptyset$ and $yz\in E(G)$, where $y\in V(I_a^\ast)$ and $z\in V(I_b^\ast)$. Assume w.l.o.g. that $y=w_2$. By Claim 1, $z\not=w_5$,  implying that either $z=w_3$ or $z=w_4$.\\

\textbf{Case 1.3.2.1}. $z=w_3$.

We have $w_1\xi_{a+1},w_1\xi_b\not\in E(G)$ (by Claim 3), $w_1w_3\not\in E(G)$ (by Claim 1) and  $w_1w_4,w_1w_5\not\in E(G)$ (by Claim 4). Assume that $N(w_1)\not\subseteq V(C)$, that is  $x_2w_1\in E(G)$ for some $x_2\in V(G\backslash C)$. Clearly, $x_2\not=x_1$, $N(x_2)\subseteq V(C)$ and $x_2w_2, x_2\xi_a\not\in E(G)$. Then we have
$$
N(x_2)\subseteq \{\xi_1,...,\xi_s,w_1\}\backslash \{\xi_a,\xi_{a+1},\xi_b\},
$$
a contradiction. Now let $N(w_1)\subseteq V(C)$. Then  
$$
N(w_1)\subseteq\{\xi_1,...,\xi_s,w_2\}\backslash \{\xi_{a+1},\xi_b\},
$$
implying that  $\xi_{a+1}=\xi_b$. Since $s\ge3$, we have $\xi_a\not=\xi_{b+1}$. Put  $I_{a-1}=\xi_{a-1}w_6\xi_a$. We have $w_6w_1,w_6w_2,w_6w_3,w_6w_5\not\in E(G)$ (by Claim 1) and  $w_6\xi_{a+1}\not\in E(G)$ (by Claim 3). If $N(w_6)\subseteq V(C)$ the  $N(w_6)\subseteq\{\xi_1,...,\xi_s\}\backslash \{\xi_{a+1}\}$, a contradiction. Analogous arguments can be used when $N(w_6)\not\subseteq V(C)$.\\

\textbf{Case 1.3.2.2}. $z=w_4$.

By Claim 5,  $w_3w_5\not\in E(G)$.   Further,  we have $w_1w_3,w_2w_5\not\in E(G)$ (by Claim 1), $w_1w_5\not\in E(G)$ (by Claim 4) and  $w_2w_3\not\in E(G)$ (otherwise we can argue as in Case 1.3.2.1). Then $G\backslash \{\xi_1,...,\xi_s,w_4\}$ has at least $s+2$ components, that is $\tau<1$, a contradiction. \\

\textbf{Case 1.3.3}. $|I_a|+|I_b|=6$.

Clearly $|I_a|=|I_b|=3$. Put $I_a=\xi_aw_1w_2\xi_{a+1}$ and $I_b=\xi_bw_3w_4\xi_{b+1}$. Assume w.l.o.g. that $w_2w_3\in E(G)$. We have $w_1\xi_{a+1},w_1\xi_b\not\in E(G)$ (by Claim 3), $w_1w_3\not\in E(G)$ (by Claim 1), $w_1w_4\not\in E(G)$ (by Lemma 2). Assume that $N(w_1)\not\subseteq V(C)$,  that is  $x_2w_1\in E(G)$ for some $x_2\in V(G\backslash C)$. Clearly, $x_2\not=x_1$, $N(x_2)\subseteq V(C)$ and $x_2w_2, x_2\xi_a\not\in E(G)$. Then we have
$$
N(x_2)\subseteq \{\xi_1,...,\xi_s,w_1\}\backslash \{\xi_a,\xi_{a+1},\xi_b\},
$$
a contradiction. Now let $N(w_1)\subseteq V(C)$. Then 
$$
N(w_1)\subseteq\{\xi_1,...,\xi_s,w_2\}\backslash \{\xi_{a+1},\xi_b\},
$$
implying that $\xi_{a+1}=\xi_b$. Since $s\ge3$, we have $\xi_a\not=\xi_{b+1}$. Put $I_{a-1}=\xi_{a-1}w_5\xi_a$. We have $w_5w_1,w_5w_2,w_5w_3,w_5w_4\not\in E(G)$ (by Claim 1), $w_5\xi_{a+1}\not\in E(G)$ (by Claim 3). If $N(w_5)\subseteq V(C)$ then $N(w_5)\subseteq\{\xi_1,...,\xi_s\}\backslash \{\xi_{a+1}\}$, a contradiction. Analogous arguments can be used when $N(w_5)\not\subseteq V(C)$.\\

\textbf{Case 1.4}. $|i:|I_i|\ge3|=1$.

Let $|I_1|\ge3$ and $|I_i|=2$ $(i=2,3,...,s)$. Clearly $3\le|I_1|\le6$.\\

\textbf{Case 1.4.1}. $|I_1|=6$.

Put $I_1=\xi_1w_1w_2w_3w_4w_5\xi_2$. \\

\textbf{Case 1.4.1.1}.  $\Upsilon(I_1,...,I_s)=\Upsilon(I_1,I_a)$ for some $a\in\{2,...,s\}$.

Clearly $|I_a|=2$.  Put $I_a=\xi_aw_6\xi_{a+1}$. We have $w_6w_1,w_6w_5\not\in E(G)$ (by Claim 1),  $w_6w_2,w_6w_4\in E(G)$ (by Claim 4) and  $w_1w_3,w_3w_5,w_1w_5\not\in E(G)$ (by Claim 5). Hence $G\backslash \{\xi_1,...,\xi_s,w_2,w_4\}$ has at least $s+3$ components, that is $\tau\le1$, a contradiction.\\

\textbf{Case 1.4.1.2}. $\Upsilon(I_1,...,I_s)\not=\Upsilon(I_1,I_i)$ for each $i\in\{2,...,s\}$.

It follows that $\Upsilon(I_1,I_a)\not=\emptyset$ and $\Upsilon(I_1,I_b)\not=\emptyset$ for some distinct $a,b\in \{2,...,s\}$. Assume w.l.o.g. that $\xi_1,\xi_a,\xi_b$ occur on $C$ in this order.Put $I_a=\xi_aw_6\xi_{a+1}$ and $I_b=\xi_bw_7\xi_{b+1}$. Let $y_1w_6,y_2w_7\in E(G)$, where $y_1,y_2\in V(I_1^\ast)$. By Claim 3, $\{y_1,y_2\}\cap\{w_1,w_5\}=\emptyset$. Further, by Claim 4, $y_1\not=y_2^+$ and $y_2\not=y_1^+$. So, $\{y_1,y_2\}=\{w_2,w_4\}$. By Claim 5, $w_1w_3,w_3w_5,w_1w_5\not\in E(G)$. This means that $G\backslash \{\xi_1,...,\xi_s,w_2,w_4\}$ has at least $s+3$ components, contradicting the fact that $\tau>1$.\\ 

\textbf{Case 1.4.2}. $|I_1|=5$.

Put $I_1=\xi_1w_1w_2w_3w_4\xi_2$. \\

\textbf{Case 1.4.2.1}.  $\Upsilon(I_1,...,I_s)=\Upsilon(I_1,I_a)$ for some $a\in\{2,...,s\}$.

Clearly $|I_a|=2$.  Put $I_a=\xi_aw_5\xi_{a+1}$. By Claim 3, $w_5w_1,w_4\not\in E(G)$. Assume w.l.o.g. that $w_5w_2\in E(G)$. But then, by Claim 4, $w_5w_3\not\in E(G)$, that is $\Upsilon(I_1,...,I_s)=\{w_5w_2\}$. So,  $G\backslash \{\xi_1,...,\xi_s,w_2\}$ has at least $s+1$ components, contradicting the fact that $\tau>1$.\\ 

\textbf{Case 1.4.2.2}. $\Upsilon(I_1,...,I_s)\not=\Upsilon(I_1,I_i)$ for each $i\in\{2,...,s\}$.

It follows that $\Upsilon(I_1,I_a)\not=\emptyset$ and $\Upsilon(I_1,I_b)\not=\emptyset$ for some distinct $a,b\in \{2,...,s\}$. Assume w.l.o.g. that $\xi_1,\xi_a,\xi_b$ occur on $C$ in this order. Put $I_a=\xi_aw_5\xi_{a+1}$ and $I_b=\xi_bw_6\xi_{b+1}$. Let $y_1w_5,y_2w_6\in E(G)$, where $y_1,y_2\in V(I_1^\ast)$. By Claim 3, $\{y_1,y_2\}\cap\{w_1,w_4\}=\emptyset$. Further, by Claim 4, $\{y_1,y_2\}\not=\{w_2,w_3\}$, that is either $\{y_1,y_2\}=\{w_2\}$ or $\{y_1,y_2\}=\{w_3\}$, say $\{y_1,y_2\}=\{w_2\}$. This means that $w_2$ belongs to all edges in $\Upsilon(I_1,...,I_s)$, implying that $G\backslash \{\xi_1,...,\xi_s,w_2\}$ has at least $s+1$ components, contradicting the fact that $\tau>1$.\\ 

\textbf{Case 1.4.3}. $|I_1|=4$.

Put $I_1=\xi_1w_1w_2w_3\xi_2$. It is not hard to see that $w_2$ belongs to all edges in $\Upsilon(I_1,...,I_s)$, implying that $G\backslash \{\xi_1,...,\xi_s,w_2\}$ has at least $s+1$ components, contradicting the fact that $\tau>1$.\\ 

\textbf{Case 1.4.4}. $|I_1|=3$.

By Lemma 2, $\Upsilon(I_1,...,I_s)=\emptyset$, implying that $\tau<1$, a contradiction.\\

\textbf{Case 2}. $\overline{p}=1$.

Since $\delta\ge\kappa\ge3$, we have $|N_C(x_i)|\ge\delta-\overline{p}=\delta-1\ge2$ \ $(i=1,2)$.\\

\textbf{Case 2.1.} $N_C(x_1)\not=N_C(x_2)$. 

It follows that $\max \{\sigma_1,\sigma_2\}\ge1$, where 
$$
\sigma_1=|N_C(x_1)\backslash N_C(x_2)|,   \quad \sigma_2=|N_C(x_2)\backslash N_C(x_1)|.
$$
If $\max \{\sigma_1,\sigma_2\}\ge3$ then by Lemma 1, $c\ge3\delta+2\ge2\delta+5$, contradicting (1). If $\max \{\sigma_1,\sigma_2\}=2$ then clearly $s\ge\delta+1$ and it is easy to see that there are at least $\delta$ elementary segments on $C$ of length at least 3. But then $c\ge2+3\delta\ge2\delta+5$, contradicting (1). Finally, let  $\max \{\sigma_1,\sigma_2\}=1$. This implies   $s\ge\delta$ and $|I_i|\ge3$ $(i=1,...,s)$. If  $s\ge \delta+1$ then  $c\ge3s\ge3\delta+3>2\delta+5$, again contradicting (1). Let $s=\delta$, that is  $|I_i|=3$ $(i=1,...,s)$. By Lemma 2, $\Upsilon(I_1,...,I_s)=\emptyset$, contradicting the fact that  $\tau>1$. \\

\textbf{Case 2.2}. $N_C(x_1)= N_C(x_2)$.

Clearly, $s=|N_C(x_1)|\ge\delta-\overline{p}=\delta-1$. If $s\ge\delta+1$ then $c\ge3s\ge3\delta+3>2\delta+5$, contradicting (1). Next suppose that $s=\delta$. If $\Upsilon(I_1,...,I_s)=\emptyset$ then $G\backslash \{\xi_1,...,\xi_s\}$ has at least $s+1$ components,  contradicting the fact that $\tau>1$. Otherwise $\Upsilon(I_a,I_b)\not=\emptyset$ for some distinct $a,b\in \{1,...,s\}$. By definition, there is an intermediate path $L$ between $I_a$ and $I_b$. By Lemma 2, 
$$
|I_a|+|I_b|\ge2\overline{p}+2|L|+4\ge8,
$$
implying that $c\ge 8+3(s-2)=3\delta+2\ge2\delta+5$, contradicting (1). So, $s=\delta-1$.   If $s=2$ then $G\backslash \{\xi_1,\xi_2\}$ is disconnected, contradicting the fact that $\kappa\ge3$. Thus $s\ge3$, implying that $\delta\ge4$.

The next claim can be derived from (1) and Lemma 2 easily.\\

\textbf{Claim 10}. $(1)$ $|I_i|+|I_j|\le9$ for each distinct $i,j\in \{1,...,s\}$.

$(2)$ If $|I_a|+|I_b|=9$ for some distinct $a,b\in \{1,...,s\}$ then  $|I_i|=3$ for each $i\in \{1,...,s\}\backslash \{a,b\}$. 

$(3)$ If $|I_a|=6$ for some $a\in\{1,...,s\}$ then $|I_i|=3$ for each $i\in \{1,...,s\}\backslash \{a\}$.  

$(4)$ There are at most three segments of length at least 4.

$(5)$ If $|I_a|\ge4$, $|I_b|\ge4$, $|I_f|\ge4$ for some distinct $a,b,f\in\{1,2,...,s\}$ then $|I_a|=|I_b|=|I_f|=4$.\\

The following three claims are the exact analogs of Claims 3,4,5 for $\overline{p}=1$ and can be proved by a similar way.\\

\textbf{Claim 11}. Let  $\xi_a,\xi_b,\xi_f$ be  distinct elements of $X$, occuring on $\overrightarrow{C}$ in a consecutive order.  If $\xi_a^-\xi_b^+\in E(G)$ then $w\xi_a, w\xi_b \not\in E(G)$ for each $w\in \{\xi_f^+,\xi_f^{++},\xi_f^-,\xi_f^{--}\}$. If either $\xi_a^{--}\xi_b^+\in E(G)$  or $\xi_a^{-}\xi_b^{++}\in E(G)$ then $w\xi_a,w\xi_b\not\in E(G)$ for each $w\in\{\xi_f^+,\xi_f^-\}$.\\

\textbf{Claim 12}. Let $a,b\in \{1,...,s\}$. If $\xi_a^+w\in E(G)$ for some $w\in V(\xi_b\overrightarrow{C}\xi_a^-)$ then $\xi_b^-y\not\in E(G)$ for each $y\in\{w^+,w^{++},w^-,w^{--}\}$ and $\xi_b^{--}y\not\in E(G)$ for each $y\in\{w^+,w^-\}$ .\\

\textbf{Claim 13}. Let $a,b\in \{1,...,s\}$. If $\xi_a^+w\in E(G)$ for some $w\in V(\xi_b^+\overrightarrow{C}\xi_a)$ then $\xi_b^+w^+,\xi_b^+w^{++},\xi_b^{++}w^+\not\in E(G)$.  \\

\textbf{Claim 14}. Let  $\xi_a,\xi_b,\xi_f,\xi_g$ be  distinct elements of $X$, occuring on $\overrightarrow{C}$ in a consecutive order. If $\xi_a^-\xi_f^+\in E(G)$ then $\xi_b^+\xi_g^-,\xi_b^{++}\xi_g^-,\xi_b^+\xi_g^{--}\not\in E(G)$.

\textbf{Proof}. If $\xi_b^+\xi_g^-\in E(G)$ then 
$$
\xi_a\overrightarrow{C}\xi_bx_2\xi_g\overrightarrow{C}\xi_a^-\xi_f^+\overrightarrow{C}\xi_g^-\xi_b^+\overrightarrow{C}\xi_fx_1\xi_a
$$
is longer than $C$, a contradiction. Hence  $\xi_b^+\xi_g^-\not\in E(G)$. Similarly, $\xi_b^{++}\xi_g^-$, $\xi_b^+\xi_g^{--}\not\in E(G)$. Claim 14 is proved.   \  \  \  $\Delta$          \\

 \textbf{Case 2.2.1}. $\mid i:\mid I_i\mid\ge4\mid=3$.  

Let $|I_i|\ge4$ for some $a,b,f\in\{1,...,s\}$ and $|I_i|=3$ for each 
$$
i\in\{1,...,s\}\backslash\{a,b,f\}.
$$
By Claim 10(5), $|I_a|=|I_b|=|I_f|=4$. Put 
$$
I_a=\xi_aw_1w_2w_3\xi_{a+1}, \  I_b=\xi_bw_4w_5w_6\xi_{b+1}, \  I_f=\xi_fw_7w_8w_9\xi_{f+1}.
$$
Assume w.l.o.g. that $\xi_a,\xi_b,\xi_f$ occur on $C$ in this order. By Lemma 2, 
$$
\Upsilon(I_1,...,I_s)=\Upsilon(I_a,I_b,I_f),
$$
$$
|\Upsilon(I_a,I_b)|\le1,\  |\Upsilon(I_a,I_f)|\le1,  \   |\Upsilon(I_b,I_f)|\le1,
$$
implying that $|\Upsilon(I_a,I_b,I_f)|\le3$. If  $|\Upsilon(I_a,I_b,I_f)|\le1$ then clearly $\tau\le1$, a contradiction. Let $|\Upsilon(I_a,I_b,I_f)|\ge2$. Assume w.l.o.g. that $\Upsilon(I_a,I_b)\not=\emptyset$ and $\Upsilon(I_a,I_f)\not=\emptyset$. By Claim 1, either $w_3w_4\in E(G)$ or $w_2w_5\in E(G)$ or $w_1w_6\in E(G)$.\\

\textbf{Case 2.2.1.1}. $w_3w_4\in E(G)$.

By Claim 12, $w_1w_9, w_2w_8\not\in E(G)$ (due to $w_3w_4\in E(G)$), implying that $w_3w_7\in E(G)$. If $\Upsilon(I_1,...,I_s)=\{w_3w_4,w_3w_7\}$ then $G\backslash \{\xi_1,...,\xi_s,w_3\}$ has at least $s+1$ components, contradicting the fact that $\tau >1$. Let $\Upsilon(I_b,I_f)\not=\emptyset$. By Claim 12, $w_6w_7, w_5w_8\not\in E(G)$ (due to $w_3w_4$), implying that $w_4w_9\in E(G)$ which contradicts Claim 14 due to $w_3w_7$ .\\

\textbf{Case 2.2.1.2}. $w_2w_5\in E(G)$.

If $w_1w_9\in E(G)$ then we can argue as in Case 2.2.1.1. Let $w_1w_9\not\in E(G)$. We have also $w_3w_7\not\in E(G)$ (by Claim 13).  Then, by Claim 1, $w_2w_8\in E(G)$. Next, by Claim 13,   $w_1w_3\not\in E(G)$.   Analogously, $w_4w_6,w_7w_9\not\in E(G)$. If 
$\Upsilon(I_1,...,I_s)=\{w_2w_5,w_2w_8\}$ then $G\backslash \{\xi_1,...,\xi_s,w_2,w_5,w_8\}$ has at  least $s+4$ components, contradicting the fact that $\tau>1$.  If $|\Upsilon(I_1,...,I_s)|=3$ then due to above observations, $w_5w_8\in E(G)$ and again $G\backslash \{\xi_1,...,\xi_s,w_2,w_5,w_8\}$ has at  least $s+4$ components, a contradiction.\\

\textbf{Case 2.2.1.3}. $w_1w_6\in E(G)$.

If either $w_1w_9\in E(G)$ or $w_2w_8\in E(G)$ then we can argue as in Cases 2.2.1.1-2.2.1.2. Otherwise, by Lemma 1,  $w_3w_7\in E(G)$ which  contradicts Claim 14 due to $w_1w_6\in E(G)$. \\

\textbf{Case 2.2.2}. $\mid i:\mid I_i\mid\ge4\mid=2$.  

Let $|I_i|\ge4$ for some $a,b\in\{1,...,s\}$ and $|I_i|=3$ for each $i\in\{1,...,s\}\backslash \{a,b\}$. By Claim 10(1), $8\le |I_a|+|I_b|\le 9$.\\

\textbf{Case 2.2.2.1}. $|I_a|+|I_b|=8$.

It follows that $|I_a|=|I_b|=4$. By Lemma 2, $\Upsilon(I_1,...,I_s)=\Upsilon(I_a,I_b)$ and $|\Upsilon(I_a,I_b)|=1$, implying that $\tau\le1$, a contradiction.\\

\textbf{Case 2.2.2.2}. $|I_a|+|I_b|=9$.

Assume w.l.o.g. that $|I_a|=4$ and $|I_b|=5$. Put $I_a=\xi_aw_1w_2w_3\xi_{a+1}$ and $I_b=\xi_bw_4w_5w_6w_7\xi_{b+1}$. If $\Upsilon(I_1,...,I_s)=\Upsilon(I_a,I_b)$ then by Lemma 2, $|\Upsilon(I_a,I_b)|\le2$ and the edges in $\Upsilon(I_a,I_b)$ have a common vertex. This means that $\tau\le1$, a contradiction. Let $\Upsilon(I_1,...,I_s)\not=\Upsilon(I_a,I_b)$. It follows that  $\Upsilon(I_b,I_f)\not=\emptyset$ for some $f\in\{1,...,s\}\backslash \{a,b\}$ and $|I_f|=3$. Put $I_f=\xi_fw_8w_9\xi_{f+1}$.\\

\textbf{Claim 15}. The edges in $\Upsilon(I_1,...,I_s)\backslash \Upsilon(I_a,I_b)$ have a common vertex.

\textbf{Proof}. By Lemma 2, $|\Upsilon(I_b,I_f)|=1$. If $\Upsilon(I_1,...,I_s)\backslash \Upsilon(I_a,I_b)=\Upsilon(I_b,I_f)$ then we are done.  Otherwise $\Upsilon(I_b,I_g)\not=\emptyset$ for some $g\in\{1,...,s\}\backslash \{a,b,f\}$. Put $I_g=\xi_gw_{10}w_{11}\xi_{g+1}$. Assume w.l.o.g. that $\xi_b,\xi_f,\xi_g$ occur on $C$ in this order. Let $y_1z_1\in\Upsilon(I_b,I_f)$ with $y_1\in V(I_b^\ast)$, $z_1\in V(I_f^\ast)$ and $y_2z_2\in\Upsilon(I_b,I_g)$ with $y_2\in V(I_b^\ast)$, $z_2\in V(I_g^\ast)$. By Claim 1, $y_1z_1\in\{w_6w_8,w_5w_9\}$ and $y_2z_2\in\{w_5w_{11},w_6w_{10}\}$. If  $y_1z_1=w_6w_8$ then by Claim 12, $y_2z_2\not=w_5w_{11}$. Then $y_2z_2=w_6w_{10}$ and we are done. If $y_1z_1=w_5w_9$ then by Claim 5,  $y_2z_2\not=w_6w_{10}$, implying that $y_2z_2=w_5w_{11}$ and again we are done. Claim 15 is proved.   \ \ \ $\Delta$\\

If $\Upsilon(I_a,I_b)=\emptyset$ then by Claim 15, there is a vertex $v$ which is incident to all edges in $\Upsilon(I_1,...,I_s)$, implying that  $G\backslash \{\xi_1,...,\xi_s,v\}$ has at least $s+1$ components, contradicting the fact that $\tau>1$. Let $\Upsilon(I_a,I_b)\not=\emptyset$.  Let $y_1z_1\in\Upsilon(I_a,I_b)$ with $y_1\in V(I_a^\ast)$, $z_1\in V(I_b^\ast)$ and $y_2z_2\in\Upsilon(I_b,I_f)$ with $y_2\in V(I_b^\ast)$, $z_2\in V(I_f^\ast)$.  By Claim 1, 
$$
y_1z_1\in\{w_1w_6,w_1w_7,w_2w_5,w_2w_6,w_3w_4,w_3w_5\},  \  \ y_2z_2\in\{w_6w_8,w_5w_9\}.
$$ 
Assume first that $y_2z_2=w_6w_8$. By Claim 12, $y_1z_1\not\in\{w_3w_4,w_3w_5,w_2w_5\}$. Next, by Claim 13, $y_1z_1\not=w_1w_7$. Hence, $y_1z_1\in\{w_1w_6,w_2w_6\}$, implying that $G\backslash \{\xi_1,...,\xi_s,w_6\}$ has at least $s+1$ components, contradicting the fact that $\tau>1$. Now let  $y_2z_2=w_5w_9$. By Claim 12, $y_1z_1\not=w_3w_4$. Further, by Claim 13, $y_1z_1\not\in \{w_1w_6,w_2w_6\}$ and by Claim 14, $y_1z_1\not=w_1w_7$. So, $y_1z_1\in\{w_2w_5,w_3w_5\}$, implying that $G\backslash \{\xi_1,...,\xi_s,w_5\}$ has at least $s+1$ components, contradicting the fact that $\tau>1$.\\

\textbf{Case 2.2.3}. $|i:|I_i|\ge4|=1$.

Let $|I_1|\ge4$ and $|I_i|=3$ $(i=2,3,...,s)$. If $|I_1|=4$ then by Lemma 2, $\Upsilon(I_1,...,I_s)=\emptyset$, implying that $\tau\le1$, a contradiction. Then by Claim 10(1), $5\le|I_1|\le6$. If $\Upsilon(I_1,...,I_s)=\Upsilon(I_1,I_a)$ for some $a\in\{2,...,s\}$ then by Lemma 2, there is a vertex which is incident to all vertices in $\Upsilon(I_1,...,I_s)$, implying that $\tau\le1$, a contradiction. Otherwise  $\Upsilon(I_1,I_a)\not=\emptyset$ and $\Upsilon(I_1,I_b)\not=\emptyset$ for some distinct $a,b\in \{2,...,s\}$. Clearly, $|I_a|=|I_b|=3$. Assume w..o.g. that $\xi_1,\xi_a,\xi_b$ occur on $C$ in this order. Let $y_1z_1\in \Upsilon(I_1,I_a)$ and $y_2z_2\in \Upsilon(I_1,I_b)$, where $y_1,y_2\in V(I_1^\ast)$. 
\\

\textbf{Case 2.2.3.1}. $|I_1|=6$.

Put 
$$
I_1=\xi_1w_1w_2w_3w_4w_5\xi_2, \  I_a=\xi_aw_6w_7\xi_{a+1},  \  I_b=\xi_bw_8w_9\xi_{b+1}. 
$$

\textbf{Claim 16}. Either $y_1=y_2$ or $y_1z_1=w_2w_7, y_2z_2=w_4w_8$.

\textbf{Proof}. By Claim 1, 
$$
y_1z_1\in\{w_3w_6,w_4w_6,w_2w_7,w_3w_7\},  \  \  y_2z_2\in\{w_2w_9,w_3w_9,w_3w_8,w_4w_8\}.
$$
By Claim 12, if $z_1=w_6$ and $z_2=w_9$ then   $y_1=y_2=w_3$.  By the same reason, if  either $z_1=w_7, z_2=w_9$ or $z_1=w_6, z_2=w_8$ then  again $y_1=y_2$. Thus, we have $y_1z_1=w_2w_7, y_2z_2=w_4w_8$. Claim 16 is proved.   \  \  \   $\Delta$\\

By Claim 16, either there is a vertex which is incident to all edges in  $\Upsilon(I_1,...,I_s)$, implying that $\tau\le1$, a contradiction, or $w_2$ and $w_4$ belong to all edges in $\Upsilon(I_1,...,I_s)$. We have $w_1w_5\not\in E(G)$ by Claim 13. Next, if $w_1w_3\in E(G)$ then 
$$
\xi_1x_1\xi_b\overleftarrow{C}\xi_{a+1}x_2\xi_1\overrightarrow{C}w_7w_2w_1w_3w_4w_8\overrightarrow{C}\xi_1
$$
is longer than $C$, a contradiction. Let $w_1w_3\not\in E(G)$. Analogously, $w_3w_5\not\in E(G)$. So, $\{w_1,w_3,w_5\}$ is an independent set of vertices and hence 
$$
G\backslash \{\xi_1,...,\xi_s,w_2,w_4\}
$$
 has at least $s+3$ components, contradicting the fact that $\tau>1$.\\

\textbf{Case 2.2.3.2}. $|I_1|=5$.

Put $I_1=\xi_1w_1w_2w_3w_4\xi_2$, $I_a=\xi_aw_5w_6\xi_{a+1}$ and $I_b=\xi_bw_7w_8\xi_{b+1}$. By Claim 1, $y_1z_1\in\{w_3w_5,w_2w_6\}$ and $y_2z_2\in \{w_2w_8,w_3w_7\}$. Then by Claim 12, $y_1=y_2$, that is there is a vertex which is incident to all edges in $\Upsilon(I_1,...,I_s)$, implying that $\tau\le1$, a contradiction.\\

\textbf{Case 3}. $\overline{p}\ge2$.

\textbf{Claim 17}. Let $x_1\xi_a,x_2\xi_{a+1}\in E(G)$ for some $a\in\{1,...,s\}$ and let $\overrightarrow{Q}=y\overrightarrow{Q}z$ be a path with
$$
y\in V(P),  \  z\in V(I_a^\ast),  \  V(Q)\cap V(P\cup C)=\{y,z\}.
$$
Then $|I_a|\ge\overline{p}+4$. If $x_1x_2\in E(G)$ and $y\not\in \{x_1,x_2\}$ then $|I_a|\ge\overline{p}+6$.

\textbf{Proof}. Since $C$ is extreme, we have 
$$
|\xi_a\overrightarrow{C}z|\ge|\xi_ax_1\overrightarrow{P}y\overrightarrow{Q}z|\ge|x_1\overrightarrow{P}y|+2,
$$
$$
|z\overrightarrow{C}\xi_{a+1}|\ge|z\overleftarrow{Q}y\overrightarrow{P}x_2\xi_{a+1}|\ge|y\overrightarrow{P}x_2|+2,
$$
implying that 
$$
|I_a|=|\xi_a\overrightarrow{C}z|+|z\overrightarrow{C}\xi_{a+1}|\ge|x_1\overrightarrow{P}y|+|y\overrightarrow{P}x_2|+4=\overline{p}+4. 
$$
Now let $x_1x_2\in E(G)$. Then 
$$
|\xi_a\overrightarrow{C}z|\ge|\xi_ax_1x_2\overleftarrow{P}y\overrightarrow{Q}z|\ge|y\overrightarrow{P}x_2|+3,
$$
$$
|z\overrightarrow{C}\xi_{a+1}|\ge|z\overleftarrow{Q}y\overleftarrow{P}x_1x_2\xi_{a+1}|\ge|x_1\overrightarrow{P}y|+3,
$$
implying that 
$$
|I_a|=|\xi_a\overrightarrow{C}z|+|z\overrightarrow{C}\xi_{a+1}|\ge|x_1\overrightarrow{P}y|+|y\overrightarrow{P}x_2|+6=\overline{p}+6. 
$$
Claim 17 is proved.     \  \  \   $\Delta$\\

\textbf{Case 3.1}. $\overline{p}\le\delta-3$.

It follows that $|N_C(x_i)|\ge\delta-\overline{p}\ge3$ $(i=1,2)$.  If $N_C(x_1)\not=N_C(x_2)$ then by  Lemma 1, $|C|\ge4\delta-2\overline{p}\ge2\delta+6$, contradicting (1). Hence $N_C(x_1)=N_C(x_2)$, implying that  $|I_i|\ge\overline{p}+2$ $(i=1,2,...,s)$. Clearly, $s\ge |N_C(x_1)|-(|V(P)|-1)\ge\delta-\overline{p}\ge3$. If $s\ge\delta-\overline{p}+1$ then 
$$
|C|\ge s(\overline{p}+2)\ge(\delta-\overline{p}+1)(\overline{p}+2)
$$
$$
=(\delta-\overline{p}-3)(\overline{p}-2)+4\delta-4\ge2\delta+6,
$$
again contradicting (1). Hence $s=\delta-\overline{p}$. It means that $x_1x_2\in E(G)$, that is $G[V(P)]$ is hamiltonian. By symmetric arguments, $N_C(y)=N_C(x_1)$ for each $y\in V(P)$. If $\Upsilon(I_1,I_2,...,I_s)=\emptyset$ then  $\tau\le1$, contradicting the hypothesis. Otherwise  $\Upsilon(I_a,I_b)\not=\emptyset$ for some elementary segments $I_a$ and $I_b$. By definition, there is an intermediate path $L$ between $I_a$ and $I_b$. If $|L|\ge2$ then by lemma 2,
$$
|I_a|+|I_b|\ge2\overline{p}+2|L|+4\ge2\overline{p}+8.
$$
Hence
$$
|C|=|I_a|+|I_b|+\sum_{i\in\{1,...,s\}\backslash \{a,b\}}|I_i|\ge2\overline{p}+8+(s-2)(\overline{p}+2)
$$
$$
=(\delta-\overline{p}-3)(\overline{p}-2)+4\delta-\overline{p}-2\ge 2\delta+6,
$$
contradicting (1).  Thus, $|L|=1$, i.e. $\Upsilon(I_1,I_2,...,I_s)\subseteq E(G)$. By Lemma 2, 
$$
|I_a|+|I_b|\ge2\overline{p}+2|L|+4= 2\overline{p}+6,
$$
which yields $|C|\ge2\delta+4$. If $|I_a|+|I_b|\ge2\overline{p}+7$, then  clearly $c\ge2\delta+5$,  contradicting (1). Hence,  If $|I_a|+|I_b|=2\overline{p}+6$ and $|I_i|=\overline{p}+2$ for each $i\in\{1,...,s\}\backslash \{a,b\}$. If $|I_a|=|I_b|=\overline{p}+3$ then by Lemma 2, $\Upsilon(I_1,...,I_s)=\Upsilon(I_a,I_b)$ and $|\Upsilon(I_a,I_b)|=1$. This means that $\tau\le1$, contradicting the fact that $\tau>1$. Now let $|I_a|=\overline{p}+4$ and $|I_b|=\overline{p}+2$. Put $I_a=\xi_aw_1w_2w_3$. Then by Claim 1, $w_2$ is incident to all edges in $\Upsilon(I_1,...,I_s)$, implying that $\tau\le1$, a contradiction. \\

\textbf{Case 3.2}. $\overline{p}=\delta-2$.

We have  $|N_C(x_i)|\ge\delta-\overline{p}=2$ $(i=1,2)$. \\

\textbf{Case 3.2.1}. $N_C(x_1)\not=N_C(x_2)$.

It follows that $s\ge3$. Clearly, there are at least two elementary segments on $C$ of length at least $\overline{p}+2$. If $s\ge5$ then $c\ge2(\overline{p}+2)+2(s-2)\ge2\delta+6$, contradicting (1). Thus, $3\le s\le4$.\\

\textbf{Case 3.2.1.1}. $s=3$.\\

\textbf{Claim 18}. For each pair $\xi_i,\xi_j$ $(i,j\in \{1,2,3\})$ there is a path with endvertices $\xi_i,\xi_j$ and vertex set $\{\xi_i,\xi_j\}\cup V(P)$.

\textbf{Proof}. Assume w.l.o.g. that $i=1, j=2$ and  $x_1\xi_1\in E(G)$. If  $x_2\xi_2\in E(G)$ then we are done. Let $x_2\xi_2\not\in E(G)$, implying that $x_1\xi_2\in E(G)$. Next, if  $x_2\xi_1\in E(G)$ then again we are done due to $\xi_1x_2\overleftarrow{P}x_1\xi_2$. Hence $x_2\xi_1\not\in E(G)$. But then $|N_C(x_2)|\le1$, a contradiction. Claim 18 is proved.  \  \  \  $\Delta$\\

By Claim 18,  $c\ge3(\overline{p}+2)=3\delta$. If $\delta\ge5$ then $c\ge3\delta\ge2\delta+5$, contradicting (1). Let $\delta\le4$. Recalling also that $\delta=\overline{p}+2\ge4$, we get  $\delta=4$, $\overline{p}=2$ and $|I_i|=\overline{p}+2=4$ $(i=1,2,3)$. Put $P=x_1yx_2$. By Claim 17, there is no a path $\overrightarrow{Q}=y\overrightarrow{Q}z$ with 
$$
z\in V(C)\backslash \{\xi_1,\xi_2,\xi_3\},  \  V(Q)\cap V(P\cup C)=\{y,z\}.
$$
Observing also that Lemma 2 is applicable in this special case due to Claim 18, we can state that $\Upsilon(I_1,I_2,I_3)=\emptyset$. But then $G\backslash \{\xi_1,\xi_2,\xi_3\}$ has at least four components, contradicting the fact that $\tau>1$.\\

\textbf{Case 3.2.1.2}. $s=4$.

If $N_C(x_1)\cap N_C(x_2)\not=\emptyset$ then clearly there are at least three elementary segments of length at least $\overline{p}+2$, which yields $c\ge3(\overline{p}+2)+2(s-3)>2\delta+5$, contradicting (1). Let $N_C(x_1)\cap N_C(x_2)=\emptyset$. Since there are at least  two elementary segments of length at least $\overline{p}+2$,  we have $c\ge2(\overline{p}+2)+2(s-2)\ge2\delta+4$. By (1), we can assume w.l.o.g. that $|I_1|=|I_3|=2$ and $|I_2|=|I_4|=\overline{p}+2$. By Claim 17, there is no a path $\overrightarrow{Q}=y\overrightarrow{Q}z$ such that 
$$
y\in V(P),  \  z\in V(C)\backslash \{\xi_1,\xi_2,\xi_3,\xi_4\},  \  V(Q)\cap V(P\cup C)=\{y,z\}.
$$
Next, if $\overrightarrow{L}=y\overrightarrow{L}z\in\Upsilon(I_1,I_2)$, where $y\in V(I_1^\ast)$ and $z\in V(I_2^\ast)$, then
$$
c\ge|\xi_1y\overrightarrow{L}z\overleftarrow{C}\xi_2x_1\overrightarrow{P}x_2\xi_3\overrightarrow{C}\xi_1|\ge|C|-|z\overrightarrow{C}\xi_3|+\overline{p}+2,
$$
implying that $|z\overrightarrow{C}\xi_3|\ge\overline{p}+2$. But Then $|I_2|>\overline{p}+2$, a contradiction. Hence $\Upsilon(I_1,I_2)=\emptyset$. Analogously, $\Upsilon(I_1,I_4)=\Upsilon(I_2,I_3)=\Upsilon(I_3,I_4)=\emptyset$. Further, if $\Upsilon(I_2,I_4)\not=\emptyset$ then we can argue as in proof of Lemma 2 to show that $|I_2|+|I_4|\ge2\overline{p}+6$, a contradiction. Hence, $\Upsilon(I_2,I_4)=\emptyset$. By a symmetric arguments, $\Upsilon(I_1,I_3)=\emptyset$. So, $\Upsilon(I_1,I_2,I_3,I_4)=\emptyset$, implying that $G\backslash \{\xi_1,\xi_2,\xi_3,\xi_4\}$ has at least five components, which contradicts the fact that $\tau>1$.\\

\textbf{Case 3.2.2}. $N_C(x_1)=N_C(x_2)$.

If $s\ge4$ then $c\ge s(\overline{p}+2)\ge4\delta\ge2\delta+8$, contradicting (1). Let $2\le s\le 3$.\\

\textbf{Case 3.2.2.1}. $s=3$.

If $\overline{p}\ge3$ then $\delta=\overline{p}+2\ge5$ and $c\ge3(\overline{p}+2)=3\delta\ge2\delta+5$, contradicting (1). Hence $\overline{p}=2$ and $\delta=4$. By (1), $c=12=2\delta+4$ and $|I_i|=4$ $(i=1,2,3)$. Put $P=x_1yx_2$. By Claim 17, there is no a path $\overrightarrow{Q}=y\overrightarrow{Q}z$ such that 
$$
  z\in V(C)\backslash \{\xi_1,\xi_2,\xi_3\},  \  V(Q)\cap V(P\cup C)=\{y,z\}.
$$
By Lemma 2, $\Upsilon(I_1,I_2,I_3)=\emptyset$, contradicting the fact that $\tau>1$.\\

\textbf{Case 3.2.2.2}. $s=2$.

It follows that $x_1x_2\in E(G)$. Since $\kappa\ge3$, there is a path $\overrightarrow{Q}=y\overrightarrow{Q}z$ such that 
$$
y\in V(P),  \    z\in V(C)\backslash \{\xi_1,\xi_2\},  \  V(Q)\cap V(P\cup C)=\{y,z\}.
$$
Assume w.l.o.g. that $z\in V(I_1^\ast)$. By Claim 17, $|I_1|\ge\overline{p}+6$. Observing also that $|I_2|\ge\overline{p}+2$, we get $c\ge2\overline{p}+8=2\delta+6$, contradicting (1). \\

\textbf{Case 3.3}. $\overline{p}=\delta-1$.

It follows that $|N_C(x_i)|\ge\delta-\overline{p}=1$ $(i=1,2)$. \\

\textbf{Case 3.3.1}. $|N_C(x_i)|\ge2$ $(i=1,2)$.

If $N_C(x_1)\not=N_C(x_2)$ then by Lemma 1, 
$|C|\ge2\overline{p}+8=2\delta+6$, contradicting (1). 
Hence, $N_C(x_1)=N_C(x_2)$. By the hypothesis,  $s\ge2$ and $|I_i|\ge\overline{p}+2=\delta+1$ $(i=1,...,s)$. If $s\ge3$ then $c\ge s(\overline{p}+2)\ge3(\delta+1)\ge2\delta+6$, contradicting (1). Let $s=2$. Since $\kappa\ge3$, there is a path $\overrightarrow{Q}=y\overrightarrow{Q}z$ such that 
$$
y\in V(P),  \    z\in V(C)\backslash \{\xi_1,\xi_2\},  \  V(Q)\cap V(P\cup C)=\{y,z\}.
$$
Assume w.l.o.g. that $z\in V(I_1^\ast)$. If $x_1x_2\in E(G)$ then by Claim 17, $|I_1|\ge\overline{p}+6$. This implies $c\ge2\overline{p}+8=2\delta+6$, contradicting (1). Let $x_1x_2\not\in E(G)$, which yields $x_1w,x_2w\in E(G)$ for each $w\in V(P)\backslash \{x_1,x_2\}$. If $\overline{p}\ge3$ then either $y^+\not\in\{x_1,x_2\}$ or $y^-\not\in\{x_1,x_2\}$, say $y^+\not\in\{x_1,x_2\}$. Then $x_1y^+\in E(G)$ and therefore,
$$
|z\overleftarrow{Q}y\overleftarrow{P}x_1y^+\overrightarrow{P}x_2\xi_2|\ge\overline{p}+2.
$$
Observe also that $y\not\in\{x_1,x_2\}$ and $|\xi_1\overrightarrow{C}z|\ge|\xi_1x_1\overrightarrow{P}y\overrightarrow{Q}z|\ge3$. Then $|I_1|\ge|\xi_1\overrightarrow{C}z|+|z\overrightarrow{C}\xi_2|\ge\overline{p}+2$ and hence $c=|I_1|+|I_2|\ge2\overline{p}+7=2\delta+5$, contradicting (1). Now let $\overline{p}=2$, implying that $\delta=3$ and 
$$
|\xi_1\overrightarrow{C}z|=|z\overrightarrow{C}\xi_2|=3,  \  |I_2|=4,  \  Q=yz,  \   |C|=10=2\delta+4.
$$
By arguing as in proof of Lemma 2, we can show that there are no edges connecting the interior vertices of the segments $\xi_1\overrightarrow{C}z$, $z\overrightarrow{C}\xi_2$ and $I_2$. Thus, $G\backslash \{\xi_1,\xi_2,z\}$ has at least four components, contradicting the fact that $\tau>1$.\\

\textbf{Case 3.3.2}. Either $|N_C(x_1)|=1$ or $|N_C(x_2)|=1$.

Assume w.l.o.g. that $|N_C(x_1)|=1$ and $N_C(x_1)=\{\xi_1\}$. It follows that $x_1x_2\in E(G)$.\\

\textbf{Case 3.3.2.1}. $N_C(x_1)\not=N_C(x_2)$.

Assume w.l.o.g. that $x_2\xi_2\in E(G)$. If $s\ge4$ then $c\ge2(\overline{p}+2)+2(s-2)\ge2\delta+6$, contradicting (1).  Let $s\le3$. \\

\textbf{Case 3.3.2.1.1}. $s=3$.

It follows that $c\ge2(\overline{p}+2)+2=2\delta+4$. By (1), $|I_1|=|I_3|=\delta+1$ and  $|I_2|=2$.  Put $I_2=\xi_2z\xi_3$. If $\Upsilon(I_1,I_2)\not=\emptyset$, that is $yz\in E(G)$ for some $y\in V(I_1^\ast)$, then
$$
|C|\ge|\xi_1x_1\overrightarrow{P}x_2\xi_2\overleftarrow{C}yz\overrightarrow{C}\xi_1|\ge|C|-|\xi_1\overrightarrow{C}y|+\overline{p}+2,
$$
implying that $|\xi_1\overrightarrow{C}y|\ge\overline{p}+2$. But then $|I_1|\ge\overline{p}+3=\delta+2$, a contradiction. Hence $\Upsilon(I_1,I_2)=\emptyset$. Analogously, $\Upsilon(I_2,I_3)=\emptyset$. If $\Upsilon(I_1,I_3)\not=\emptyset$ then we can argue as in proof of Lemma 2, to show that $|I_1|+|I_2|\ge2\overline{p}+6=2\delta+4$, a contradiction. So, $\Upsilon(I_1,I_2,I_3)=\emptyset$. Further, if there is a path $\overrightarrow{Q}=w_1\overrightarrow{Q}w_2$ such that 
$$
w_1\in V(P),  \    w_2\in V(C)\backslash \{\xi_1,\xi_2,\xi_3\},  \  V(Q)\cap V(P\cup C)=\{w_1,w_2\},
$$
then clearly $w_2\not\in V(I_2^\ast)$ (since $|I_2|=2$) and $w_2\not\in V(I_1^\ast)\cup V(I_3^\ast)$ by Claim 17. Otherwise, $G\backslash \{\xi_1,\xi_2,\xi_3\}$ has at least four components, contradicting the fact that $\tau>1$.\\

\textbf{Case 3.3.2.1.2}. $s=2$.

Since $\kappa\ge3$, there is a path $\overrightarrow{Q}=y\overrightarrow{Q}z$ such that 
$$
y\in V(P),  \    z\in V(C)\backslash \{\xi_1,\xi_2\},  \  V(Q)\cap V(P\cup C)=\{y,z\}.
$$
Clearly, $y\not\in\{x_1,x_2\}$. Assume w.l.o.g. that $z\in V(I_1^\ast)$. By Claim 17, $|I_1|\ge\overline{p}+6$, implying that $c\ge2\overline{p}+8=2\delta+6$, contradicting (1).\\

 \textbf{Case 3.3.2.2}. $N_C(x_1)=N_C(x_2)$.

It follows that $N_C(x_1)=N_C(x_2)=\{x_1\}$ and $x_1w\in E(G)$ for each $w\in V(P)\backslash \{x_1\}$. Since $\kappa\ge3$, there is a path $\overrightarrow{Q}=y\overrightarrow{Q}z$ such that 
$$
y\in V(P),  \    z\in V(C)\backslash \{\xi_1\},  \  V(Q)\cap V(P\cup C)=\{y,z\}.
$$
Clearly, $y\not\in \{x_1,x_2\}$. Since $x_1y^+\in E(G)$, we can replace $P$ with $y\overleftarrow{P}x_1y^+\overrightarrow{P}x_2$. Then we can argue as in Case 3.3.2.1.\\

\textbf{Case 3.4}. $\overline{p}\ge\delta$.

If $|C|\ge\kappa(\delta+1)$ then clearly $|C|\ge3(\delta+1)\ge2\delta+6$, contradicting (1). Otherwise, by Lemma 3,  $|N_C(x_i)|\ge2$ $(i=1,2)$. Since there are at least two elementary segments on $C$ of length at least $\overline{p}+2\ge\delta+2$, we have   $|C|\ge2(\overline{p}+2)+2\ge2\delta+6$ when $s\ge3$, contradicting (1). Now let $s=2$. By (1), 
$$
\overline{p}=\delta,  \  |I_1|=|I_2|=\delta+2,  \  c=2\delta+4.
$$
Since $\kappa\ge3$, there is a path $\overrightarrow{Q}=y\overrightarrow{Q}z$ such that 
$$
y\in V(P),  \    z\in V(C)\backslash \{\xi_1,\xi_2\},  \  V(Q)\cap V(P\cup C)=\{y,z\}.
$$
Assume w.l.o.g. that $z\in V(I_1^\ast)$. But then, by Claim 17, $|I_1|\ge\overline{p}+4=\delta+4$, a contradiction.            \quad   \quad         \rule{7pt}{6pt}

\noindent Institute for Informatics and Automation Problems\\ National Academy of Sciences\\
P. Sevak 1, Yerevan 0014, Armenia\\ E-mail: zhora@ipia.sci.am
\end{document}